\documentclass{amsart}
\usepackage{amsfonts}
\usepackage{amssymb}

\newtheorem{theorem}{Theorem}[section]

\newtheorem{conjecture}[theorem]{Conjecture}
\newtheorem{corollary}[theorem]{Corollary}

\newtheorem{definition}[theorem]{Definition}

\newtheorem{lemma}[theorem]{Lemma}

\newtheorem{proposition}[theorem]{Proposition}

\numberwithin{equation}{section}

\begin{document}
\title{Poisson Equation on complete manifolds}
\author{Ovidiu Munteanu}
\email{ovidiu.munteanu@uconn.edu}
\address{Department of Mathematics, University of Connecticut, Storrs, CT
06268, USA}
\author{Chiung-Jue Anna Sung}
\email{cjsung@math.nthu.edu.tw}
\address{Department of Mathematics, National Tsing Hua University, Hsin-Chu,
Taiwan}
\author{Jiaping Wang}
\email{jiaping@math.umn.edu}
\address{School of Mathematics, University of Minnesota, Minneapolis, MN
55455, USA}
\thanks{The first author was partially supported by NSF grant DMS-1506220.
The second author was partially supported by MOST. The third author was
partially supported by NSF grant DMS-1606820 }

\begin{abstract}
We develop heat kernel and Green's function estimates for manifolds with
positive bottom spectrum. The results are then used to establish existence
and sharp estimates of the solution to the Poisson equation on such
manifolds with Ricci curvature bounded below. As an application, we show
that the curvature of a steady gradient Ricci soliton must decay
exponentially if it decays faster than linear and the potential function is
bounded above.
\end{abstract}

\maketitle

\section{Introduction}

For a complete noncompact manifold $(M^{n},g)$ without boundary, consider
the Poisson equation

\begin{equation*}
\Delta u=-\varphi ,
\end{equation*}%
where $\varphi$ is a given smooth function on $M.$ In this paper, we
establish existence and sharp estimates of the solution $u$ and provide
applications to steady gradient Ricci solitons.

As well-known, the solvability of the Poisson equation is closely related to
the existence of the so-called Green's function. In \cite{M}, Malgrange
showed that $M$ always admits a Green's function $G(x,y),$ namely, $%
G(x,y)=G(y,x)$ and $\Delta _{y}G(x,y)=-\delta _{x}(y).$ In particular, if $%
\varphi \in C_{0}^{\infty }(M),$ then a solution $u$ to the Poisson equation
exists and is given by

\begin{equation*}
u(x)=\int_{M}G(x,y)\,\varphi (y)\,dy.
\end{equation*}

Malgrange's proof is rather abstract. Later, Li and Tam \cite{LT} provided a
more constructive proof. Among other things, the constructed Green's
function satisfies

\begin{equation*}
\sup_{y\in M\setminus B(p,2R)}\sup_{x\in B(p,R)}|G(p,y)-G(x,y)|<\infty .
\end{equation*}%
This turns out to be very useful in applications. For example, it was used
to prove an extension theorem for harmonic functions in \cite{STW}.

\begin{theorem}
(Sung-Tam-Wang) For any harmonic function $u$ defined on $M\setminus \Omega, 
$ where $\Omega$ is a bounded subset of $M,$ there exists a harmonic
function $v$ on $M$ such that $u-v$ is bounded on $M\setminus \Omega.$
\end{theorem}

Obviously, a good control of the Green's function $G(x,y)$ will enable one
to establish existence and estimates of the solution $u$ to the Poisson
equation for more general $\varphi.$ Recall that $M$ is nonparabolic if $M$
admits a positive Green's function and parabolic otherwise. It is well-known
that $M$ is nonparabolic if and only if it admits a nonconstant bounded
superharmonic function. Therefore, when $M$ is parabolic, the solution $u$
to the Poisson equation must be unbounded if $\varphi \geq 0$ but not identically $0.$ Since we are
primarily concerned on the existence of bounded solutions in this paper, we
will restrict our attention to nonparabolic manifold $M.$ In particular,
there exists a unique minimal positive Green's function on $M.$ In the
following, it is understood that this is the Green's function we refer to.

For manifolds with nonnegative Ricci curvature, a sharp pointwise estimate
for the Green's function is available.

\begin{theorem}
(Li-Yau) Let $M^n$ be a complete manifold with nonnegative Ricci curvature.
If $\int_{1}^{\infty }\mathrm{V}^{-1}(p,\sqrt{t})\,dt<\infty $ for some
point $p\in M,$ then $M$ is nonparabolic and its minimal positive Green's
function $G(x,y)$ satisfies the estimate

\begin{equation*}
C_{1}\,\int_{r^{2}(x,y)}^{\infty }\mathrm{V}^{-1}(x,\sqrt{t})dt\leq
G(x,y)\leq C_{2}\,\int_{r^{2}(x,y)}^{\infty }\mathrm{V}^{-1}(x,\sqrt{t})dt
\end{equation*}
for some constants $C_{1}$ and $C_{2}$ depending only on the dimension $n.$
\end{theorem}

Here and in the following, $\mathrm{V}(p,r)$ denotes the volume of the
geodesic ball $B(p,r)$ centered at point $p$ with radius $r,$ and $r(x,y)$ the distance
between points $x$ and $y$ in $M.$ The estimate
follows from their famous upper and lower bounds of the heat kernel \cite{LY}
together with the fact that

\begin{equation*}
G(x,y)=\int_0^{\infty} H(x,y,t) dt.
\end{equation*}

Based on the estimates of the Green's function, Ni, Shi and Tam \cite{NST}
obtained the following result concerning the Poisson equation.

\begin{theorem}
(Ni-Shi-Tam) Let $M$ be a complete manifold with nonnegative Ricci
curvature. If $M$ is nonparabolic, then for a locally H\"{o}lder continuous
function $\varphi $ with $|\varphi |(x)\leq c\,r^{-k}(x)$ for some $k>2,$
the Poisson equation $\Delta u=-\varphi $ has a solution $u$ such that $%
|u|(x)\leq C\,r^{-k+2}(x).$
\end{theorem}

In fact, they have proved more general results. For the existence, the decay
rate on $\varphi $ is only assumed to be $k>1.$ Moreover, the
decay on $\varphi $ is only required to be true in the average sense over
the geodesic balls centered at a fixed point. While the solution $u$ in
general is no longer bounded, its growth is well controlled.

They applied their result to study, among other things, the following
uniformization conjecture of Yau.

\begin{conjecture}
(Yau) A complete noncompact K\"ahler manifold with positive bisectional curvature is
biholomorphic to the complex Euclidean space.
\end{conjecture}

Indeed, by first solving the Poisson equation $\Delta u=S$ on such manifold $%
M,$ where $S$ is the scalar curvature, they demonstrated that under suitable
assumptions $u$ is in fact a solution to the Poincar\'e-Lelong equation

\begin{equation*}
\sqrt{-1}\partial \bar{\partial}u=\rho,
\end{equation*}
where $\rho$ is the Ricci form of $M.$

This line of ideas was initiated by Mok, Siu and Yau in \cite{MSY}. While
the conjecture in its most general form is still open, there are various
partial results. We refer to the recent spectacular work of Liu \cite{L} and
the references therein for further information.

Our focus will be on manifolds with positive spectrum. Denote by $%
\lambda_{1}(\Delta )$ the smallest spectrum of the Laplacian or the bottom
spectrum of $M.$ It is well-known that $M$ is nonparabolic if $%
\lambda_{1}(\Delta )>0.$ Recall that $\lambda _{1}(\Delta )$ can be
characterized as the best constant of the Poincar\'{e} inequality.

\begin{equation*}
\lambda _{1}(\Delta )=\inf_{\phi \in C_{0}^{\infty }}\frac{\int_{M}|\nabla
\phi |^{2}dx}{\int_{M}\phi ^{2}dx}.
\end{equation*}

As observed by Strichartz \cite{S}, if $\lambda _{1}(\Delta )>0,$ then $%
\Delta ^{-1}$ is in fact a bounded operator on $L^{p}(M)$ for $1<p<\infty.$
In particular, there exists a solution $u\in L^{p}(M)$ to the Poisson
equation for $\varphi \in L^{p}(M).$

Our achievement here is to establish an existence result with sharp control
of the solution $u$ by only assuming a modest decay on the function $%
\varphi, $ very much in the spirit of the result alluded above by Ni, Shi
and Tam for the case of nonnegative Ricci curvature.

\begin{theorem}
Let $M$ be a complete Riemannian manifold with bottom spectrum $\lambda
_{1}(\Delta )>0$ and Ricci curvature $\mathrm{Ric}\geq -\left( n-1\right) K$
for some constant $K.$ Let $\varphi $ be a smooth function such that

\begin{equation*}
\left\vert \varphi \right\vert \left( x\right) \leq c\,\left(
1+r(x)\right)^{-k}
\end{equation*}
for some $k>1,$ where $r(x)$ is the distance function from $x$ to a fixed
point $p\in M.$ Then the Poisson equation $\Delta u=-\varphi$ admits a
bounded solution $u$ on $M.$ If, in addition, the volume of the ball $B(x,1)$
satisfies $\mathrm{V}\left( x,1\right) \geq c$ for all $x\in M,$ then the
solution $u$ decays and

\begin{equation*}
\left\vert u\right\vert \left( x\right) \leq C\,\left( 1+r(x)\right) ^{-k+1}.
\end{equation*}
\end{theorem}

We point out that the existence of a solution $u$ was previously proved by
the first author and Sesum \cite{MS}. However, their estimate on the
solution $u$ takes the form

\begin{equation*}
\left\vert u\right\vert \left( x\right) \leq C\,e^{c\,r(x)}.
\end{equation*}

It should also be emphasized that the assumption on the volume that $\mathrm{%
V}(x,1)\geq c$ is necessary to guarantee the solution $u$ decays at
infinity. Indeed, since $u$ is a bounded super-harmonic function when $%
\varphi $ is positive, $u$ can not possibly decay to $0$ along a parabolic
end of $M.$

The theorem is sharp as one can see as follows. On the hyperbolic space $%
\mathbb{H}^{n},$ the Green's function is given by

\begin{equation*}
G\left( x,y\right) =\int_{r\left( x,y\right) }^{\infty }\frac{dt}{A\left(
t\right) },
\end{equation*}
where $A\left( t\right) $ is the area of geodesic sphere of radius $t$ in $%
\mathbb{H}^{n}.$ For $\varphi \left( x\right) =\left( 1+r(x)\right) ^{-k}$
with $k>1,$ a direct calculation gives

\begin{eqnarray*}
u\left( x\right) &=&\int_{\mathbb{H}^{n}}G\left( x,y\right) \varphi \left(
y\right) dy \\
&\geq &c\,\left( 1+r(x)\right) ^{-k+1}.
\end{eqnarray*}

Our proof again relies on some sharp estimates of the Green's function.
Recall the following result of the third author with Li \cite{LW}, which is
a sharp version of Agmon's work \cite{A}.

\begin{theorem}
(Li-Wang) Let $M$ be a complete Riemannian manifold with $\lambda
_{1}(\Delta )>0.$ Let $u$ be a nonnegative subharmonic function defined on $%
M\setminus \Omega ,$ where $\Omega $ is a compact domain. If $u$ satisfies
the growth condition

\begin{equation*}
\int_{(M\setminus \Omega )\cap B(p,R)}u^{2}\,e^{-2\sqrt{\lambda _{1}(\Delta )%
}\,r}=o(R)
\end{equation*}
as $R\rightarrow \infty,$ then it must satisfy the decay estimate

\begin{equation*}
\int_{B(p,R+1)\setminus B(p,R)}u^{2}\leq C\,e^{-2\sqrt{\lambda _{1}(\Delta )}%
\,R}
\end{equation*}
for some constant $C>0$ depending on $u$ and $\lambda _{1}(\Delta ).$
\end{theorem}

In particular, the theorem implies that the minimal positive Green's
function satisfies

\begin{equation*}
\int_{B(p,R+1)\setminus B(p,R)}G^{2}(p,y)\,dy\leq C\,e^{-2\sqrt{\lambda
_{1}(\Delta )}\,R}.
\end{equation*}

While this result provides a version of sharp estimate on the Green's
function, to prove our theorem, however, we also need the following double
integral estimate.

\begin{equation*}
\int_{A}\int_{B}G(x,y)\,dy\,dx \leq \frac{e^{\sqrt{\lambda _{1}\left(
\Delta\right) }}}{\lambda _{1}\left( \Delta \right) }\,\sqrt{\mathrm{V}%
\left( A\right) }\sqrt{ \mathrm{V}\left( B\right) }\left( 1+\,r(A,B)\right)
e^{-\sqrt{\lambda _{1}\left( \Delta \right) }\,r(A,B)}
\end{equation*}
for any bounded domains $A$ and $B$ of $M.$

For this purpose, we develop a parabolic version of the aforementioned
result of Li and the third author.

\begin{theorem}
Let $M$ be a complete Riemannian manifold with $\lambda _{1}(\Delta )>0.$
Suppose that $(\Delta -\frac{\partial }{\partial t})u(x,t)\geq 0$ with $%
u(x,t)\geq 0,$

\begin{equation*}
\int_{M}u^{2}(x,0)e^{2\sqrt{\lambda _{1}(\Delta)}\,r(x,A)}dx<\infty
\end{equation*}%
and 
\begin{equation*}
\int_{0}^{T}\int_{B(A, 2R)\setminus B(A,R)}u^{2}(x,t)\,e^{-2\sqrt{\lambda
_{1}(\Delta)}\,r(x,A)}dxdt=o(R)
\end{equation*}%
for all $T>0$ as $R\rightarrow \infty.$ Then, for all $R>0,$

\begin{equation*}
\int_{0}^{\infty }\int_{B(A,R+2)\setminus B(A,R)}u^{2}(x,t)\,dxdt\leq
C\,e^{-2\sqrt{\lambda _{1}(\Delta)}\,R}\int_{M}u^{2}(x,0)e^{2\sqrt{\lambda
_{1}(\Delta)}\,r(x,A)}dx.
\end{equation*}
\end{theorem}

Here, $A$ is a bounded subset in $M,$ $r(x,A)$ the distance from $x$ to $A$
and $B(A,R)=\{x\in M\,|\,r(x,A)<R\}.$ The constant $C>0$ depends only on $%
\lambda _{1}\left( \Delta \right).$

By applying the theorem to the function $u(x,t)=\int_{A}H(x,y,t)dy,$
one obtains a sharp integral estimate of the heat kernel. This may be of
independent interest. The desired estimate of the Green's function 
follows from the fact that

\begin{equation*}
G(x,y)=\int_0^{\infty} H(x,y,t)\,dt.
\end{equation*}

The following result is crucial to our proof of Theorem 1.5. It
provides a sharp integral control of the Green's function.

\begin{theorem}
Let $M^n$ be an $n$-dimensional complete manifold with $\lambda _{1}\left(
\Delta \right) >0$ and $\mathrm{Ric}\geq -\left( n-1\right) K.$ Then for any 
$x\in M$ and $r>0$ we have

\begin{equation*}
\int_{B\left( p,r\right) }G\left( x,y\right) dy\leq C\,\left( 1+r\right)
\end{equation*}%
for some constant $C$ depending on $n,$ $K$ and $\lambda _{1}\left( \Delta
\right).$ 
\end{theorem}

On top of the double integral estimate, the proof of the theorem utilizes an
idea originated in \cite{LW1} and further illuminated in \cite{MS}, where
they used the co-area formula together with suitably chosen cut-off
functions to justify that for any $x\in M$ and $0<\alpha <\beta,$

\begin{equation*}
\int_{L_{x}\left( \alpha ,\beta \right) }G\left( x,y\right) dy\leq c\,\left(
1+\ln \frac{\beta }{\alpha }\right).
\end{equation*}%
Here, 
\begin{equation*}
L_{x}\left( \alpha ,\beta \right) :=\left\{ y\in M:\alpha <G\left(
x,y\right) <\beta \right\} .
\end{equation*}

Partly motivated by applications to gradient Ricci solitons, we in fact consider 
more generally the weighted Poisson equation on
smooth metric measure space $(M,g,e^{-f}\,dx),$ that is, Riemannian manifold 
$(M,g)$ together with a weighted measure $e^{-f}\,dx,$ where $f$ is a smooth
function on $M.$ The weighted Poisson equation is given by 

\begin{equation*}
\Delta _{f}u=-\varphi,
\end{equation*}%
where $\Delta _{f}\,u=\Delta u-\left\langle \nabla f,\nabla u\right\rangle $
is the weighted Laplacian. A natural curvature notion corresponding to the
Ricci curvature in the Riemannian setting is the Bakry-Emery Ricci
curvature, defined by

\begin{equation*}
\mathrm{Ric}_{f}=\mathrm{Ric}+\mathrm{Hess}\left( f\right).
\end{equation*}

It is known (see \cite{MW1}) that results such as volume comparison,
gradient estimates and mean value inequality are available on the smooth
metric measure spaces with the Bakry-Emery Ricci curvature bounded below
together with suitable assumptions on the weight function $f.$ With this in
mind, we have a parallel version of Theorem 1.5.

\begin{theorem}
\label{intro}Let $\left( M^n,\,g,\,e^{-f}dx\right) $ be an $n$-dimensional
smooth metric measure space with $\mathrm{Ric}_{f}\geq -\left( n-1\right) K$
and the oscillation of $f$ on any unit ball $B\left( x,1\right) $ bounded
above by a fixed constant $a.$ Assume that the bottom spectrum of the
weighted Laplacian $\lambda _{1}\left( \Delta _{f}\right) $ is positive. Let 
$\varphi $ be a smooth function such that

\begin{equation*}
\left\vert \varphi \right\vert \left( x\right) \leq c\,\left( 1+r(x)\right)
^{-k}
\end{equation*}
for some $k>1.$ Then $\Delta _{f}u=-\varphi $ admits a bounded solution $u$
on $M.$ If, in addition, the weighted volume of the ball $B(x,1)$ satisfies $%
\mathrm{V}_{f}\left( x,1\right) \geq c$ for all $x\in M,$ then the solution $%
u$ decays and satisfies

\begin{equation*}
\left\vert u\right\vert \left( x\right) \leq C\,\left( 1+r(x)\right) ^{-k+1}.
\end{equation*}
\end{theorem}

In fact, we establish a slightly more general result (see Theorem \ref{Poisson1}).
As an immediate application, we obtain the following decay estimate concerning
the subsolutions to semi-linear equations. It would be interesting to see if the estimate
can be improved to exponential decay.

\begin{theorem} Let $\left( M^n,\,g,\,e^{-f}dx\right) $ be an $n$-dimensional
smooth metric measure space with $\mathrm{Ric}_{f}\geq -\left( n-1\right) K$
and the oscillation of $f$ on any unit ball $B\left( x,1\right) $ bounded
above by a fixed constant $a.$ Assume that the bottom spectrum $\lambda _{1}\left( \Delta _{f}\right) $
of the weighted Laplacian is positive and
the weighted volume has lower bound $\mathrm{V}_{f}\left( x,1\right) \geq c>0$
for all $x\in M.$ Suppose $\psi \geq 0$ satisfies

\begin{equation*}
\Delta _{f}\psi \geq -c\psi ^{q}
\end{equation*}%
for some $q>1,$ and 
\begin{equation*}
\lim_{x\rightarrow \infty }\psi \left( x\right) r^{\frac{1}{q-1}}\left(
x\right) =0.
\end{equation*}%
Then there exist $\delta >0$ and $C>0$ such that 
\begin{equation*}
\psi \left( x\right) \leq Ce^{-r^{\delta }\left( x\right) }.
\end{equation*}
\end{theorem}

This result motivated us to study the curvature behavior of steady gradient
Ricci solitons. 

\begin{definition}
A steady gradient Ricci soliton is a complete manifold $\left( M,g\right)$
on which there exists a smooth potential function $f$ such that 
\begin{equation*}
\mathrm{Ric}+\mathrm{Hess}\left( f\right) =0.
\end{equation*}
\end{definition}

Steady gradient Ricci solitons are self-similar solutions to the Ricci flow.
Indeed, if we let $g(t)=\psi(t)^{*}\,g,$ where $\psi(t)$ is the
diffeomorphism generated by the vector field $\nabla f$ with $\psi(0)=id_M,$
then $g(t)$ is a solution to the Ricci flow 
\begin{equation*}
\frac{\partial}{\partial t}\,g(t)=-2 \mathrm{Ric}(g(t)).
\end{equation*}
As such, they play important role in the study of the Ricci flows. Some
prominent examples of steady gradient Ricci solitons include the Euclidean
space $\mathbb{R}^{n}$ with $f$ being a linear function, Hamilton's cigar
soliton $(\Sigma,\,g_{\Sigma }),$ where $\Sigma =\mathbb{R}^{2}$ and

\begin{equation*}
g_{\Sigma }=\frac{dx^{2}+dy^{2}}{1+x^{2}+y^{2}}
\end{equation*}
with the potential function $f(x,y)=-\ln (1+x^{2}+y^{2}),$ and Bryant
soliton $(\mathbb{R}^{n},g),$ $n\geq 3,$ where $g$ is rotationally symmetric
and $f=f(r)$ as well. The scalar curvature of the cigar satisfies $S=e^{f}$
and decays exponentially $S\simeq ce^{-r\left( x\right) }$ in the distance
function. However, the curvature of the Bryant soliton decays linearly in
distance.

For a steady gradient Ricci soliton, its Riemann curvature $\mathrm{Rm}$  
satisfies 

\begin{equation*}
\Delta _{f}|\mathrm{Rm}| \geq -c\,|\mathrm{Rm}|^{2}
\end{equation*}
for some constant $c>0.$ Moreover, $|\nabla f|$ is bounded and $\lambda _{1}(\Delta _{f})>0$ 
by \cite{MW2}. So Theorem 1.10 becomes applicable once the weighted volume assumption 
is verified. This more or less follows from potential $f$ being bounded above by a constant.
These considerations motivate the following theorem. 

\begin{theorem}
Let $\left( M^{n},g,f\right) $ be a complete steady gradient Ricci soliton
with potential $f$ bounded above by a constant. If its Riemann curvature satisfies

\begin{equation*}
\left\vert \mathrm{Rm}\right\vert \left( x\right) r\left( x\right) =o\left(
1\right) 
\end{equation*}%
as $x\to \infty,$ then

\begin{equation*}
\left\vert \mathrm{Rm}\right\vert \left( x\right) \leq c\,\left(
1+r(x)\right) ^{3\left( n+1\right) }\,e^{-r(x)}.
\end{equation*}
\end{theorem}

It is unclear at this point whether the assumption on $f$ is necessary. It is known that the
assumption automatically holds true when $\mathrm{Ric}>0.$ We also note that the exponential 
decay rate in the theorem is sharp as seen
from $M=N\times \Sigma,$ where $\Sigma $ is the cigar soliton and $N$ a
compact Ricci flat manifold.

It view of our result, one may wonder whether there is a dichotomy for the
curvature decay rate of steady gradient Ricci solitons, namely, either
exactly linear or exponential. This dichotomy, if confirmed, should be very
useful for the classification of steady gradient Ricci solitons. In the
three dimensional case, very recently, Deng and Zhu \cite{DZ} have shown
that such a soliton must be the Bryant soliton if its curvature decays
exactly linearly. On the other hand, if the curvature decays faster than
linear, then it must be the product of the cigar soliton and a
circle (see Corollary 5.5). We should also mention that Brendle \cite{B}
has confirmed Perelman's assertion in \cite{P} that a noncollapsed three
dimensional steady gradient Ricci soliton must be the Bryant soliton.

\textbf{Acknowledgement:} We would like to dedicate this paper to Professor
Peter Li on the occasion of his sixty-fifth birthday. It can not be
overstated how much we have benefited from his teaching, encouragement and
support over the years.

\section{Heat kernel estimates}

In this section we extend the decay estimate for subharmonic functions
developed in \cite{LW} and \cite{LW1} to the subsolutions of the heat
equation. As a consequence, we obtain heat kernel estimate on complete
manifolds with positive bottom spectrum. The estimate will be applied in
next section to derive integral estimates for the minimal Green's function.

We will cast our result in a more general setting of smooth metric measure
space $\left( M,g,e^{-f}dx\right) ,$ where the following weighted Poincar%
\'{e} inequality holds true for a positive function $\rho.$

\begin{equation}
\int_{M}\rho (x)\,\phi ^{2}(x)e^{-f}\,dx\leq \int_{M}|\nabla \phi
|^{2}(x)e^{-f}\,dx  \label{a1}
\end{equation}
for any compactly supported function $\phi \in C_{0}^{\infty }(M).$

Let us define the $\rho $-metric by

\begin{equation*}
ds_{\rho }^{2}=\rho \,ds^{2}.
\end{equation*}%
Using this metric, we consider the $\rho $-distance function defined to be

\begin{equation*}
r_{\rho }(x,y)=\inf_{\gamma }\ell _{\rho }(\gamma ),
\end{equation*}%
the infimum of the length of all smooth curves joining $x$ and $y$ with
respect to $ds_{\rho }^{2}.$ For a fixed point $p\in M,$ one checks readily
that $|\nabla r_{\rho }|^{2}(p,x)=\rho (x).$ We say that manifold $M$ has
property $(P_{\rho })$ if the $\rho $-metric is complete, and this will be
our standing assumption in this section.

Similarly, for a compact domain $A\subset M,$ we denote

\begin{equation*}
r_{\rho }(x,A)=\inf_{y\in A}r_{\rho }(y,x)
\end{equation*}%
to be the $\rho $-distance to $A$ and 
\begin{equation*}
B_{\rho }(A,R)=\{x\in M\,|\,r_{\rho }(x,A)<R\}
\end{equation*}%
to be the set of points in $M$ that have $\rho $-distance less than $R$ from
set $A.$

Consider $u\left( x,t\right)$ a nonnegative subsolution to the weighted heat
equation

\begin{equation}
\left(\Delta _{f}-\frac{\partial}{\partial t}\right) u\geq 0.  \label{a2}
\end{equation}%
We assume that $u\left( x,t\right) $ satisfies the growth condition that

\begin{equation}
\int_{M}u^{2}\left( x,0\right) e^{2r_{\rho }\left( x,A\right) }e^{-f\left(
x\right) }dx<\infty  \label{a3'}
\end{equation}%
and that for all $T>0,$

\begin{equation}
\int_{0}^{T}\int_{B_{\rho}(A, 2R)\setminus B_{\rho}(A, R)}\rho \left(
x\right) u^{2}\left( x,t\right) e^{-2r_{\rho }\left( x,A\right) }e^{-f\left(
x\right) }dxdt=o(R)  \label{a3''}
\end{equation}
as $R\rightarrow \infty.$

\begin{theorem}
\label{HK}Let $\left( M,g,e^{-f}dx\right) $ be a complete smooth metric
measure space with property ($P_{\rho }$). Let $u(x,t)$ satisfy (\ref{a2}), (%
\ref{a3'}) and (\ref{a3''}). Then for all $R>0,$

\begin{eqnarray*}
&&\int_{0}^{\infty }\int_{B_{\rho }(A,R+2)\setminus B_{\rho }(A,R)}\rho
\left( x\right) \,u^{2}\left( x,t\right) \,e^{-f\left( x\right) }dx\,dt \\
&\leq &C\,e^{-2R}\int_{M}u^{2}\left( x,0\right) e^{2r_{\rho }\left(
x,A\right) }e^{-f\left( x\right) }\,dx
\end{eqnarray*}%
for some absolute constant $C>0.$
\end{theorem}

\begin{proof}
Throughout the proof, we will denote by $C$ an absolute constant which may
change from line to line. We also suppress the dependency of $A$ and write $%
B_{\rho }(R)=B_{\rho }(A,R)$ and $r_{\rho }(x)=r_{\rho }(x,A).$ The first
step is to prove that for any $0<\delta <1,$ there exists a constant $%
0<C<\infty $ such that

\begin{equation*}
\int_{0}^{\infty }\int_{M}\rho (x)e^{2\delta r_{\rho
}(x)}u^{2}(x,t)e^{-f\left( x\right) }\,dx\,dt\leq \frac{C\,}{1-\delta }%
\int_{M}u^{2}\left( x,0\right) e^{2r_{\rho }\left( x\right) }e^{-f\left(
x\right) }\,dx.
\end{equation*}

Indeed, let $\phi (x)$ be a non-negative cut-off function on $M.$ Then for
any function $h(x)$ integration by parts yields%
\begin{eqnarray}
\int_{M}|\nabla (\phi \,u\,e^{h})|^{2}e^{-f} &=&\int_{M}|\nabla (\phi
e^{h})|^{2}\,u^{2}e^{-f}+\int_{M}(\phi \,e^{h})^{2}\,|\nabla u|^{2}e^{-f}
\label{a4} \\
&&+2\int_{M}\left( \phi \,e^{h}\right) u\langle \nabla (\phi \,e^{h}),\nabla
u\rangle e^{-f}  \notag \\
&=&\int_{M}|\nabla (\phi \,e^{h})|^{2}u^{2}e^{-f}+\int_{M}\phi ^{2}|\nabla
u|^{2}e^{2h}e^{-f}  \notag \\
&&+\frac{1}{2}\int_{M}\langle \nabla (\phi ^{2}\,e^{2h}),\nabla u^{2}\rangle
e^{-f}  \notag \\
&=&\int_{M}|\nabla (\phi e^{h})|^{2}u^{2}e^{-f}+\int_{M}\phi ^{2}|\nabla
u|^{2}e^{2h}e^{-f}  \notag \\
&&-\frac{1}{2}\int_{M}\phi ^{2}\Delta _{f}(u^{2})e^{2h}e^{-f}  \notag \\
&=&\int_{M}|\nabla (\phi \,e^{h})|^{2}\,u^{2}e^{-f}-\int_{M}\phi
^{2}\,u\left( \Delta _{f}u\right) e^{2h}e^{-f}  \notag \\
&\leq &\int_{M}|\nabla \phi |^{2}\,u^{2}e^{2h}e^{-f}+2\int_{M}\phi \,\langle
\nabla \phi ,\nabla h\rangle u^{2}e^{2h}e^{-f}  \notag \\
&&+\int_{M}\phi ^{2}\,|\nabla h|^{2}\,u^{2}\,e^{2h}e^{-f}-\int_{M}\phi
^{2}uu_{t}e^{2h}e^{-f},  \notag
\end{eqnarray}%
where in the last line we have used (\ref{a2}). On the other hand, using the
weighted Poincar\'e inequality (\ref{a1}), we have

\begin{equation*}
\int_{M}\rho \,\phi ^{2}\,u^{2}e^{2h}e^{-f}\leq \int_{M}|\nabla (\phi
\,ue^{h})|^{2}e^{-f}.
\end{equation*}%
Hence (\ref{a4}) becomes

\begin{eqnarray}
\int_{M}\rho \,\phi ^{2}\,u^{2}\,e^{2h}e^{-f} &\leq &\int_{M}|\nabla \phi
|^{2}\,u^{2}\,e^{2h}e^{-f}+2\int_{M}\phi \,\langle \nabla \phi ,\nabla
h\rangle u^{2}e^{2h}e^{-f}  \label{a5} \\
&&+\int_{M}\phi ^{2}\,|\nabla h|^{2}\,u^{2}e^{2h}e^{-f}-\frac{1}{2}\frac{d}{%
dt}\int_{M}\phi ^{2}u^{2}e^{2h}e^{-f}.  \notag
\end{eqnarray}
Integrating with respect to $t,$ we conclude

\begin{eqnarray}
&&\int_{0}^{T}\int_{M}\rho \,\phi ^{2}\,u^{2}\,e^{2h}\,e^{-f}dx\,dt+\frac{1}{%
2}\int_{M}\phi ^{2}\,u^{2}(x,T)e^{2h}\,e^{-f}dx  \label{a6} \\
&\leq &\int_{0}^{T}\int_{M}|\nabla \phi
|^{2}\,u^{2}\,e^{2h}e^{-f}dx\,dt+2\int_{0}^{T}\int_{M}\phi \,u^{2}\langle
\nabla \phi ,\nabla h\rangle e^{2h}\,e^{-f}\,dx\,dt  \notag \\
&&\quad +\int_{0}^{T}\int_{M}\phi ^{2}\,|\nabla
h|^{2}u^{2}\,e^{2h}e^{-f}dx\,dt+\frac{1}{2}\,\int_{M}\phi ^{2}u^{2}\left(
x,0\right) e^{2h}e^{-f}\,dx.  \notag
\end{eqnarray}
Let us first choose

\begin{equation}
\phi \left( r_{\rho }\left( x\right) \right) =\left\{ 
\begin{array}{c}
1 \\ 
R^{-1}(2R-r_{\rho }(x)) \\ 
0%
\end{array}%
\right. 
\begin{array}{c}
\text{on} \\ 
\text{on} \\ 
\text{on}%
\end{array}%
\begin{array}{l}
B_{\rho }\left( R\right) \\ 
B_{\rho }(2R)\backslash B_{\rho }(R) \\ 
M\backslash B_{\rho }(2R)%
\end{array}
\label{phi}
\end{equation}%
and%
\begin{equation*}
h\left( r_{\rho }\left( x\right) \right) =\left\{ 
\begin{array}{c}
\delta r_{\rho }\left( x\right) \\ 
K-r_{\rho }\left( x\right)%
\end{array}%
\right. 
\begin{array}{l}
\text{on }B_{\rho }\left( \left( 1+\delta \right) ^{-1}K\right) \\ 
\text{on }M\backslash B_{\rho }\left( \left( 1+\delta \right) ^{-1}K\right)%
\end{array}%
\end{equation*}%
for some fixed $K>1$. Note that when $R\geq \left( 1+\delta \right) ^{-1}K,$ 
\begin{equation*}
\left\vert \nabla \phi \right\vert ^{2}\left( x\right) =\left\{ 
\begin{array}{c}
R^{-2}\rho \left( x\right) \\ 
0%
\end{array}%
\right. 
\begin{array}{c}
\text{on} \\ 
\text{on}%
\end{array}%
\begin{array}{l}
B_{\rho }(2R)\backslash B_{\rho }(R) \\ 
\left( M\backslash B_{\rho }(2R)\right) \cup B_{\rho }\left( R\right)%
\end{array}%
\end{equation*}%
and 
\begin{equation*}
\left\langle \nabla \phi ,\nabla h\right\rangle \left( x\right) =\left\{ 
\begin{array}{c}
R^{-1}\rho \left( x\right) \\ 
0%
\end{array}%
\right. 
\begin{array}{c}
\text{on} \\ 
\text{on}%
\end{array}%
\begin{array}{l}
B_{\rho }(2R)\backslash B_{\rho }(R) \\ 
\left( M\backslash B_{\rho }(2R)\right) \cup B_{\rho }\left( R\right) ,%
\end{array}%
\end{equation*}%
whereas%
\begin{equation*}
\left\vert \nabla h\right\vert ^{2}\left( x\right) =\left\{ 
\begin{array}{c}
\delta ^{2}\rho \\ 
\rho%
\end{array}%
\right. 
\begin{array}{l}
\text{on }B_{\rho }\left( \left( 1+\delta \right) ^{-1}K\right) \\ 
\text{on }M\backslash B_{\rho }\left( \left( 1+\delta \right) ^{-1}K\right) .%
\end{array}%
\end{equation*}
Substituting all these into (\ref{a6}) implies

\begin{eqnarray*}
&&\int_{0}^{T}\int_{M}\rho \,\phi ^{2}\,u^{2}e^{2h}e^{-f}dx\,dt+\frac{1}{2}%
\,\int_{M}\phi ^{2}u^{2}(x,T)e^{2h}\,e^{-f}dx \\
&\leq &R^{-2}\int_{0}^{T}\int_{B_{\rho }(2R)\backslash B_{\rho }(R)}\rho
\,u^{2}e^{2h}e^{-f}\,dx\,dt \\
&&+2R^{-1}\int_{0}^{T}\int_{B_{\rho }(2R)\backslash B_{\rho }(R)}\rho
\,u^{2}e^{2h}e^{-f}\,dx\,dt \\
&&+\delta ^{2}\int_{0}^{T}\int_{B_{\rho }((1+\delta )^{-1}K)}\rho \,\phi
^{2}\,u^{2}\,e^{2h}e^{-f}\,dx\,dt \\
&&+\int_{0}^{T}\int_{B_{\rho }(2R)\backslash B_{\rho }((1+\delta
)^{-1}K)}\rho \,\phi ^{2}\,u^{2}\,e^{2h}e^{-f}\,dx\,dt \\
&&+\frac{1}{2}\,\int_{M}u^{2}\left( x,0\right) e^{2h}e^{-f}\,dx.
\end{eqnarray*}
This proves that

\begin{eqnarray*}
&&(1-\delta ^{2})\int_{0}^{T}\int_{B_{\rho }((1+\delta )^{-1}K)}\rho
\,u^{2}e^{2h}e^{-f}\,dx\,dt \\
&&+\frac{1}{2}\,\int_{B_{\rho }((1+\delta )^{-1}K)}u^{2}(x,T)e^{2h}e^{-f}\,dx
\\
&\leq &R^{-2}\int_{0}^{T}\int_{B_{\rho }(2R)\backslash B_{\rho }(R)}\rho
\,u^{2}e^{2h}\,e^{-f}dx\,dt \\
&&+2R^{-1}\int_{0}^{T}\int_{B_{\rho }(2R)\backslash B_{\rho }(R)}\rho
\,u^{2}e^{2h}\,e^{-f}dx\,dt \\
&&+\frac{1}{2}\,\int_{M}u^{2}\left( x,0\right) e^{2r_{\rho }}e^{-f}\,dx.
\end{eqnarray*}%
In view of the definition of $h$ and (\ref{a3''}), the first two terms on
the right hand side of this inequality tend to $0$ as $R\rightarrow \infty.$
Therefore, we obtain the estimate

\begin{eqnarray*}
&&(1-\delta ^{2})\int_{0}^{T}\int_{B_{\rho }((1+\delta )^{-1}K)}\rho
\,u^{2}e^{2\delta r_{\rho }}e^{-f}dx\,dt \\
&&+\frac{1}{2}\,\int_{B_{\rho }((1+\delta )^{-1}K)}u^{2}(x,T)e^{2\delta
r_{\rho }}e^{-f}\,dx \\
&\leq &\frac{1}{2}\,\,\int_{M}u^{2}\left( x,0\right) e^{2r_{\rho
}}e^{-f}\,dx.
\end{eqnarray*}
Since the right hand side is independent of $K$, by letting $K\rightarrow
\infty$ we conclude that

\begin{equation}
\int_{0}^{T}\int_{M}\rho \,u^{2}e^{2\delta r_{\rho }}e^{-f}\,dx\,dt\leq 
\frac{1}{2(1-\delta ^{2})}\int_{M}u^{2}\left( x,0\right) e^{2r_{\rho
}}e^{-f}\,dx  \label{a7}
\end{equation}%
and

\begin{equation*}
\int_{M}u^{2}(x,T)\,e^{2\delta r_{\rho }}e^{-f}dx\leq \int_{M}u^{2}\left(
x,0\right) e^{2r_{\rho }}e^{-f}\,dx
\end{equation*}%
for all $T>0$ and $0<\delta <1.$

Our next step is to improve this estimate by setting $h=r_{\rho }$ in the
preceding argument. Note that with this choice of $h,$ (\ref{a5}) asserts
that

\begin{eqnarray*}
-2\,\int_{M}\phi \langle \nabla \phi ,\nabla r_{\rho }\rangle
\,u^{2}e^{2r_{\rho }}e^{-f} &\leq &\int_{M}|\nabla \phi
|^{2}u^{2}\,e^{2r_{\rho }}e^{-f} \\
&&-\frac{1}{2}\frac{d}{dt}\int_{M}\phi ^{2}u^{2}e^{2r_{\rho }}e^{-f}.
\end{eqnarray*}%
For $0<R_{1}<R$, let us choose $\phi $ to be

\begin{equation*}
\phi \left( x\right) =\left\{ 
\begin{array}{c}
R_{1}^{-1}r_{\rho }(x) \\ 
\left( R-R_{1}\right) ^{-1}\left( R-r_{\rho }(x)\right)%
\end{array}%
\right. 
\begin{array}{c}
\text{on} \\ 
\text{on}%
\end{array}%
\begin{array}{l}
B_{\rho }(R_{1})\newline
\\ 
B_{\rho }(R)\backslash B_{\rho }(R_{1}).%
\end{array}%
\end{equation*}%
We conclude that

\begin{eqnarray*}
&&\frac{2}{R-R_{1}}\,\int_{B_{\rho }(R)\backslash B_{\rho }(R_{1})}\frac{%
R-r_{\rho }}{R-R_{1}}\,\rho u^{2}\,e^{2r_{\rho }}e^{-f} \\
&\leq &\frac{2}{R_{1}^{2}}\int_{B_{\rho }(R_{1})}r_{\rho }\,\rho
\,u^{2}\,e^{2r_{\rho }}e^{-f}+\frac{1}{(R-R_{1})^{2}}\,\int_{B_{\rho
}(R)\backslash B_{\rho }(R_{1})}\rho \,u^{2}\,e^{2r_{\rho }}e^{-f} \\
&&+\frac{1}{R_{1}^{2}}\int_{B_{\rho }(R_{1})}\rho \,u^{2}\,e^{2r_{\rho
}}e^{-f}-\frac{1}{2}\frac{d}{dt}\int_{M}\phi ^{2}u^{2}e^{2r_{\rho }}e^{-f}.
\end{eqnarray*}%
Integrating with respect to $t,$ we obtain%
\begin{eqnarray*}
&&\frac{2}{R-R_{1}}\,\int_{0}^{T}\int_{B_{\rho }(R)\backslash B_{\rho
}(R_{1})}\frac{R-r_{\rho }}{R-R_{1}}\,\rho \,u^{2}\,e^{2r_{\rho
}}e^{-f}\,dx\,dt \\
&\leq &\frac{2}{R_{1}^{2}}\int_{0}^{T}\int_{B_{\rho }(R_{1})}r_{\rho }\,\rho
\,u^{2}\,e^{2r_{\rho }}e^{-f}dx\,dt \\
&&+\frac{1}{(R-R_{1})^{2}}\,\int_{0}^{T}\int_{B_{\rho }(R)\backslash B_{\rho
}(R_{1})}\rho \,u^{2}\,e^{2r_{\rho }}e^{-f}\,dx\,dt \\
&&+\frac{1}{R_{1}^{2}}\int_{0}^{T}\int_{B_{\rho }(R_{1})}\rho
\,u^{2}e^{2r_{\rho }}e^{-f}\,dx\,dt+\frac{1}{2}\int_{M}u^{2}\left(
x,0\right) e^{2r_{\rho }}e^{-f}\,dx.
\end{eqnarray*}
On the other hand, for any $0<\tau <R-R_{1},$ since

\begin{eqnarray*}
&&\frac{2\tau }{(R-R_{1})^{2}}\,\int_{B_{\rho }(R-\tau )\backslash B_{\rho
}(R_{1})}\rho \,u^{2}\,e^{2r_{\rho }}e^{-f} \\
&\leq &\frac{2}{(R-R_{1})^{2}}\,\int_{B_{\rho }(R)\backslash B_{\rho
}(R_{1})}(R-r_{\rho })\,\rho \,u^{2}e^{2r_{\rho }}e^{-f},
\end{eqnarray*}%
we deduce that

\begin{eqnarray}
&&\frac{2\tau }{(R-R_{1})^{2}}\,\int_{0}^{T}\int_{B_{\rho }(R-\tau
)\backslash B_{\rho }(R_{1})}\rho \,u^{2}e^{2r_{\rho }}e^{-f}\,dx\,dt
\label{a8} \\
&\leq &\left( \frac{2}{R_{1}}+\frac{1}{R_{1}^{2}}\right)
\int_{0}^{T}\int_{B_{\rho }(R_{1})}\rho \,u^{2}e^{2r_{\rho }}e^{-f}\,dx\,dt 
\notag \\
&&+\frac{1}{(R-R_{1})^{2}}\int_{0}^{T}\int_{B_{\rho }(R)\backslash B_{\rho
}(R_{1})}\rho \,u^{2}e^{2r_{\rho }}\,e^{-f}dx\,dt  \notag \\
&&+\frac{1}{2}\,\int_{M}u^{2}\left( x,0\right) e^{2r_{\rho }}e^{-f}\,dx. 
\notag
\end{eqnarray}%
Take $R_{1}=1,$ $\tau =1$, and set

\begin{equation*}
g(R)=\int_{0}^{T}\int_{B_{\rho }(R)\backslash B_{\rho }(1)}\rho
\,u^{2}e^{2r_{\rho }}e^{-f}\,dx\,dt.
\end{equation*}%
Using (\ref{a7}) for $\delta =0$ we may rewrite the inequality (\ref{a8}) as

\begin{equation*}
g(R-1)\leq C_{1}\,R^{2}\int_{M}u^{2}\left( x,0\right) e^{2r_{\rho
}}e^{-f}\,dx+\frac{1}{2}\,g(R),
\end{equation*}%
where 
\begin{equation*}
C_{1}=\frac{3}{4}\,\left( e^{2}+1\right)
\end{equation*}%
is an absolute constant. Iterating this inequality, we obtain that for any
positive integer $k$ and $R\geq 1,$

\begin{eqnarray*}
g(R) &\leq &C_{1}\,\left( \sum_{i=1}^{k}\frac{(R+i)^{2}}{2^{i-1}}\right)
\int_{M}u^{2}\left( x,0\right) e^{2r_{\rho }}e^{-f}\,dx+2^{-k}\,g(R+k) \\
&\leq &C\,R^{2}\int_{M}u^{2}\left( x,0\right) e^{2r_{\rho
}}e^{-f}\,dx+2^{-k}\,g(R+k)
\end{eqnarray*}%
for some absolute constant $C.$ However, our previous estimate (\ref{a7})
asserts that

\begin{equation*}
\int_{0}^{T}\int_{M}\rho \,u^{2}e^{2\delta r_{\rho }}e^{-f}\,dx\,dt\leq 
\frac{1}{2\left( 1-\delta ^{2}\right) }\int_{M}\,u^{2}\left( x,0\right)
e^{2r_{\rho }}e^{-f}\,dx
\end{equation*}%
for any $\delta <1.$ This implies that

\begin{eqnarray*}
g(R+k) &=&\int_{0}^{T}\int_{B_{\rho }(R+k)\backslash B_{\rho }(1)}\rho
\,u^{2}e^{2r_{\rho }}e^{-f}dx\,dt \\
&\leq &e^{2(R+k)(1-\delta )}\int_{0}^{T}\int_{B_{\rho }(R+k)\backslash
B_{\rho }(1)}\rho \,u^{2}e^{2\delta r_{\rho }}e^{-f}\,dx\,dt \\
&\leq &\frac{1}{2\left( 1-\delta ^{2}\right) }e^{2(R+k)(1-\delta
)}\int_{M}\,u^{2}\left( x,0\right) e^{2r_{\rho }}e^{-f}\,dx.
\end{eqnarray*}%
Hence, 
\begin{equation*}
2^{-k}\,g(R+k)\rightarrow 0
\end{equation*}%
as $k\rightarrow \infty $ by choosing $2(1-\delta )<\ln 2.$ This proves the
estimate%
\begin{equation*}
g(R)\leq C\,\,R^{2}\int_{M}u^{2}\left( x,0\right) e^{2r_{\rho }}e^{-f}\,dx.
\end{equation*}%
By adjusting the constant, we have%
\begin{equation}
\int_{0}^{T}\int_{B_{\rho }(R)}\rho \,u^{2}\,e^{2r_{\rho
}}e^{-f}\,dx\,dt\leq C\,R^{2}\int_{M}u^{2}\left( x,0\right) e^{2r_{\rho
}}e^{-f}\,dx  \label{a9}
\end{equation}%
for all $R>1.$

Using inequality (\ref{a8}) again and choosing $R_{1}=1$ and $\tau =\frac{R}{%
2}$ this time, we conclude that%
\begin{eqnarray*}
R\,\int_{0}^{T}\int_{B_{\rho }(\frac{R}{2})\backslash B_{\rho }(1)}\rho
\,u^{2}e^{2r_{\rho }}e^{-f}\,dx\,dt &\leq &CR^{2}\int_{M}u^{2}\left(
x,0\right) e^{2r_{\rho }}e^{-f}\,dx \\
&&+\int_{0}^{T}\int_{B_{\rho }(R)\backslash B_{\rho }(1)}\rho
\,u^{2}e^{2r_{\rho }}e^{-f}dx\,dt.
\end{eqnarray*}%
However, applying the estimate (\ref{a9}) to the second term on the right
hand side, we have 
\begin{equation*}
\int_{0}^{T}\int_{B_{\rho }(\frac{R}{2})\backslash B_{\rho }(1)}\rho
\,u^{2}e^{2r_{\rho }}e^{-f}dx\,dt\leq CR\int_{M}u^{2}\left( x,0\right)
e^{2r_{\rho }}e^{-f}\,dx.
\end{equation*}%
Therefore, for $R>1,$ 
\begin{equation}
\int_{0}^{T}\int_{B_{\rho }(R)}\rho \,u^{2}e^{2r_{\rho }}e^{-f}\,dx\,dt\leq
C\,R\int_{M}u^{2}\left( x,0\right) e^{2r_{\rho }}e^{-f}\,dx.  \label{a10}
\end{equation}

We are now ready to prove the theorem by using (\ref{a10}). Setting $\tau =2$
and $R_{1}=R-4$ in (\ref{a8}), we obtain%
\begin{eqnarray*}
&&\int_{0}^{T}\int_{B_{\rho }(R-2)\backslash B_{\rho }(R-4)}\rho
\,u^{2}e^{2r_{\rho }}e^{-f}dx\,dt \\
&\leq &\left( \frac{8}{R-4}+\frac{4}{(R-4)^{2}}\right)
\,\int_{0}^{T}\int_{B_{\rho }(R-4)}\rho \,u^{2}e^{2r_{\rho }}e^{-f}\,dx\,dt
\\
&&+\frac{1}{4}\int_{0}^{T}\int_{B_{\rho }(R)\backslash B_{\rho }(R-4)}\rho
\,u^{2}e^{2r_{\rho }}e^{-f}dx\,dt \\
&&+2\,\int_{M}u^{2}\left( x,0\right) e^{2r_{\rho }}e^{-f}\,dx.
\end{eqnarray*}%
According to (\ref{a10}), the first term of the right hand side is bounded by%
\begin{eqnarray*}
&&\left( \frac{8}{R-4}+\frac{4}{(R-4)^{2}}\right)
\,\int_{0}^{T}\int_{B_{\rho }(R-4)}\rho \,u^{2}e^{2r_{\rho }}e^{-f}\,dx\,dt
\\
&\leq &C\,\int_{A}u^{2}\left( x,0\right) e^{2r_{\rho }}e^{-f}\,dx.
\end{eqnarray*}%
Hence, the above inequality can be rewritten as%
\begin{eqnarray*}
\int_{0}^{T}\int_{B_{\rho }(R-2)\backslash B_{\rho }(R-4)}\rho
\,u^{2}e^{2r_{\rho }}e^{-f}dx\,dt &\leq &\frac{1}{3}\int_{0}^{T}\int_{B_{%
\rho }(R)\backslash B_{\rho }(R-2)}\rho \,u^{2}e^{2r_{\rho }}e^{-f}\,dx\,dt
\\
&&+C\,\int_{M}u^{2}\left( x,0\right) e^{2r_{\rho }}e^{-f}\,dx.
\end{eqnarray*}%
Iterating this inequality $k$ times, we arrive at 
\begin{eqnarray*}
&&\int_{0}^{T}\int_{B_{\rho }(R+2)\backslash B_{\rho }(R)}\rho
\,u^{2}e^{2r_{\rho }}e^{-f}\,dx\,dt \\
&\leq &3^{-k}\int_{0}^{T}\int_{B_{\rho }(R+2(k+1))\backslash B_{\rho
}(R+2k)}\rho \,u^{2}e^{2r_{\rho }}e^{-f}\,dx\,dt \\
&&+\,C\left( \sum_{i=0}^{k-1}3^{-i}\right) \int_{M}u^{2}\left( x,0\right)
e^{2r_{\rho }}e^{-f}\,dx\,.
\end{eqnarray*}%
However, using (\ref{a10}) again, we conclude that the second term is
bounded by 
\begin{eqnarray*}
&&3^{-k}\int_{0}^{T}\int_{B_{\rho }(R+2(k+1))\backslash B_{\rho }(R+2k)}\rho
\,u^{2}e^{2r_{\rho }}dx\,dt \\
&\leq &C\,3^{-k}(R+2(k+1))\int_{M}u^{2}\left( x,0\right) e^{2r_{\rho
}}e^{-f}\,dx
\end{eqnarray*}%
and tends to $0$ as $k\rightarrow \infty.$ Hence,

\begin{equation}
\int_{0}^{T}\int_{B_{\rho }(R+2)\backslash B_{\rho }(R)}\rho
\,u^{2}e^{2r_{\rho }}e^{-f}\,dx\,dt\leq C\,\int_{M}u^{2}\left( x,0\right)
e^{2r_{\rho }}e^{-f}\,dx  \label{a11}
\end{equation}%
for some absolute constant $C>0.$ The theorem now follows from (\ref{a11})
by letting $T\rightarrow \infty.$
\end{proof}

We now apply this theorem to establish integral estimates for the heat
kernel. For a compact set $A\subset M,$ let

\begin{equation}
u(x,t)=\int_{A}H(x,y,t)\,e^{-f\left( y\right) }dy,  \label{a12}
\end{equation}%
where $H(x,y,t)$ is the minimal heat kernel of $\Delta _{f}$ on $M.$
Clearly, $u(x,0)=\chi_{A}(x).$ Furthermore, as

\begin{equation*}
\frac{d}{dt}\int_{M}u^{2}\left( x,t\right) e^{-f\left( x\right)
}dx=-\int_{M}\left\vert \nabla u\right\vert ^{2}\left( x,t\right)
e^{-f\left( x\right) }dx\leq 0,
\end{equation*}%
it follows that%
\begin{eqnarray}
\int_{M}u^{2}(x,t)e^{-f\left( x\right) }\,dx &\leq
&\int_{M}u^{2}(x,0)e^{-f\left( x\right) }\,dx  \label{a13} \\
&\leq &\mathrm{V}_{f}(A),  \notag
\end{eqnarray}%
where $\mathrm{V}_{f}(A):=\int_{A}e^{-f}dv$ is the weighted volume of the
set $A.$ Theorem \ref{HK} then implies the following.

\begin{corollary}
\label{HKE}Let $\left( M,g,e^{-f}dx\right) $ be a complete smooth metric
measure space with property ($P_{\rho }$). Then $u\left( x,t\right) $
defined in (\ref{a12}) satisfies 
\begin{equation*}
\int_{0}^{\infty }\int_{B_{\rho }\left( A,R+1\right) \backslash B_{\rho
}\left( A,R\right) }\rho \left( x\right) u^{2}\left( x,t\right) e^{-f\left(
x\right) }dxdt\leq C\,e^{-2R}\,\mathrm{V}_{f}\left( A\right)
\end{equation*}%
for all $R>0,$ where $C>0$ is an absolute constant. Furthermore, for all $%
0<\delta<1,$

\begin{equation*}
\int_{0}^{\infty }\int_{M}\rho \left( x\right) \,u^{2}\left( x,t\right)
e^{2\delta r_{\rho }\left( x,A\right) }e^{-f\left( x\right) }\,dx\,dt\leq 
\frac{1}{2(1-\delta ^{2})}\mathrm{V}_{f}\left( A\right).
\end{equation*}
\end{corollary}

\begin{proof}
To apply Theorem \ref{HK} we need to verify the assumptions (\ref{a3'}) and (%
\ref{a3''}). Since $u\left( x,0\right) =\chi _{A}\left( x\right) $, it
follows that 
\begin{equation}
\int_{M}u^{2}\left( x,0\right) e^{2r_{\rho }\left( x,A\right) }e^{-f\left(
x\right) }dx=\mathrm{V}_{f}\left( A\right) .  \label{a14}
\end{equation}%
Choosing $h=0$ and $\phi $ as in (\ref{phi}), we get from (\ref{a6}) that 
\begin{eqnarray*}
\int_{0}^{T}\int_{B_{\rho }\left( R\right) }\rho u^{2}e^{-f}dxdt &\leq &%
\frac{1}{R^{2}}\int_{0}^{T}\int_{B_{\rho }\left( 2R\right) \backslash
B_{\rho }\left( R\right) }u^{2}e^{-f}dxdt \\
&&+\frac{1}{2}\int_{B_{\rho }\left( 2R\right) }u^{2}\left( x,0\right)
e^{-f}dx \\
&\leq &\left( \frac{T}{R^{2}}+\frac{1}{2}\right) \mathrm{V}_{f}(A),
\end{eqnarray*}%
where we have used (\ref{a13}) in the last line. Letting $R\rightarrow
\infty,$ one sees that 
\begin{equation*}
\int_{0}^{T}\int_{M}\rho u^{2}e^{-2r_{\rho }}e^{-f}dxdt<\infty
\end{equation*}
and (\ref{a3''}) follows.

Now the first estimate of the corollary follows from (\ref{a11}) and (\ref%
{a14}), and the second from (\ref{a7}) and (\ref{a14}).
\end{proof}

Finally, in the case $\lambda _{1}\left( \Delta _{f}\right)>0,$ obviously
one may take $\rho=\lambda _{1}\left( \Delta _{f}\right).$

\begin{corollary}
\label{lambda}Let $\left( M,g,e^{-f}dx\right) $ be a smooth metric measure
space with $\lambda _{1}\left( \Delta _{f}\right) >0.$ Then the function $%
u\left( x,t\right) $ defined in (\ref{a12}) satisfies

\begin{equation*}
\int_{0}^{\infty }\int_{B\left( A,R+1\right) \backslash B\left( A,R\right)
}u^{2}\left( x,t\right) e^{-f\left( x\right) }dxdt\leq C\,e^{-2\sqrt{\lambda
_{1}\left( \Delta _{f}\right) }R}\,\mathrm{V}_{f}\left( A\right),
\end{equation*}%
where $C$ depends only on $\lambda _{1}\left( \Delta _{f}\right).$
Furthermore, for all $0<\delta<1,$

\begin{equation*}
\int_{0}^{\infty }\int_{M}\,u^{2}\left( x,t\right) e^{2\delta \sqrt{\lambda
_{1}\left( \Delta _{f}\right) }r\left( x,A\right) }e^{-f\left( x\right)
}\,dx\,dt\leq \frac{1}{2(1-\delta ^{2})\lambda _{1}\left( \Delta _{f}\right) 
}\,\mathrm{V}_{f}\left( A\right).
\end{equation*}
\end{corollary}

\begin{proof}
Since $ds_{\rho }^{2}=\lambda _{1}\left( \Delta _{f}\right) ds^{2},$ we get $%
r_{\rho }\left( x,y\right) =\sqrt{\lambda _{1}\left( \Delta _{f}\right) }\,
r\left( x,y\right)$ and $B_{\rho }\left( A,R\right) =B\left( A,\frac{R}{%
\sqrt{\lambda _{1}\left( \Delta _{f}\right) }}\right).$ The result now
follows from Corollary \ref{HKE}. The constant $C$ can be taken as $C=c\cdot
\max \left\{ 1,\lambda _{1}\left( \Delta _{f}\right) ^{-1}\right\},$ where $%
c $ is an absolute constant.
\end{proof}

\section{Green's function estimates}

In this section, we develop estimates for the Green's function. The results,
while of independent interest, will be applied to solve the Poisson equation
in next section. We continue to assume that $\left( M,g,e^{-f}dx\right)$ is
a smooth metric measure space with positive bottom spectrum $\lambda
_{1}\left( \Delta _{f}\right)>0.$ It is known (see e.g. Chapter 20 in \cite%
{Li}) that $M$ must be $f$-nonparabolic, that is, there exists a positive
Green's function for the weighted Laplacian $\Delta _{f}.$ Denote by $%
G\left( x,y\right)$ the minimal positive Green's function. We have the
following.

\begin{theorem}
\label{GAB}Let $\left( M,g,e^{-f}dx\right) $ be a complete smooth metric
measure space with positive bottom spectrum $\lambda _{1}\left( \Delta
_{f}\right) .$ Then the minimal positive Green's function $G(x,y)$ of $%
\Delta _{f}$ satisfies

\begin{eqnarray*}
&&\int_{A}\int_{B}G(x,y)\,e^{-f\left( x\right) }e^{-f\left( y\right) }dy\,dx
\\
&\leq &\frac{e^{\sqrt{\lambda _{1}\left( \Delta _{f}\right) }}}{\lambda
_{1}\left( \Delta _{f}\right) }\,\sqrt{\mathrm{V}_{f}\left( A\right) }\,%
\sqrt{\mathrm{V}_{f}\left( B\right) }\,\left( 1+\,r(A,B)\right) \,e^{-\sqrt{%
\lambda _{1}\left( \Delta _{f}\right) }\,r(A,B)}
\end{eqnarray*}%
for any bounded domains $A$ and $B$ of $M.$
\end{theorem}

\begin{proof}
By Corollary \ref{lambda}, for any $0<\delta <1,$

\begin{eqnarray*}
&&\int_{0}^{\infty }\int_{M}\left( \int_{A}H(x,y,t)\,e^{-f\left( y\right)
}dy\right) ^{2}e^{2\delta \,\sqrt{\lambda _{1}\left( \Delta _{f}\right) }%
\,r(x,A)}e^{-f\left( x\right) }\,dx\,dt \\
&\leq &\frac{1}{2(1-\delta ^{2})\lambda _{1}\left( \Delta _{f}\right) }\,%
\mathrm{V}_{f}(A).
\end{eqnarray*}%
Of course, the same inequality holds for domain $B$ as well. Therefore,
noting that $r(A,B)\leq r(x,A)+r(x,B),$ we get

\begin{eqnarray*}
&&\int_{0}^{\infty }\int_{M}\left( \int_{A}H(x,y,t)e^{-f\left( y\right)
}\,dy\right) \, \\
&&\times \left( \int_{B}H(x,z,t)\,e^{-f\left( z\right) }dz\right)
\,e^{\delta \sqrt{\lambda _{1}\left( \Delta _{f}\right) }\,r(A,B)}e^{-f%
\left( x\right) }\,dx\,dt \\
&\leq &\left( \int_{0}^{\infty }\int_{M}\left( \int_{A}H(x,y,t)e^{-f\left(
y\right) }\,dy\right) ^{2}\,e^{2\delta \,\sqrt{\lambda _{1}\left( \Delta
_{f}\right) }r(x,A)}e^{-f\left( x\right) }\,dx\,dt\right) ^{1/2} \\
&&\times \left( \int_{0}^{\infty }\int_{M}\left(
\int_{B}H(x,z,t)\,e^{-f\left( z\right) }dz\right) ^{2}\,e^{2\delta \,\sqrt{%
\lambda _{1}\left( \Delta _{f}\right) }\,r(x,B)}e^{-f\left( x\right)
}\,dx\,dt\right) ^{1/2} \\
&\leq &\frac{1}{2(1-\delta ^{2})\lambda _{1}\left( \Delta _{f}\right) }\,%
\sqrt{\mathrm{V}_{f}(A)}\sqrt{\,\mathrm{V}_{f}(B)}.
\end{eqnarray*}%
However, 
\begin{eqnarray*}
&&\int_{0}^{\infty }\int_{M}\left( \int_{A}H(x,y,t)e^{-f\left( y\right)
}dy\right) \\
&&\times \,\left( \int_{B}H(x,z,t)\,e^{-f\left( z\right) }dz\right)
\,e^{\delta \sqrt{\lambda _{1}\left( \Delta _{f}\right) }\,r(A,B)}e^{-f%
\left( x\right) }\,dx\,dt \\
&=&\int_{A}\int_{B}\int_{0}^{\infty }\int_{M}H(x,y,t)\,H(x,z,t)e^{-f\left(
x\right) }e^{-f\left( y\right) }e^{-f\left( z\right) }\,e^{\delta \sqrt{%
\lambda _{1}\left( \Delta _{f}\right) }\,r(A,B)}dx\,dt\,dz\,dy\, \\
&=&\int_{A}\int_{B}\int_{0}^{\infty }H(y,z,2t)\,e^{-f\left( y\right)
}e^{-f\left( z\right) }e^{\delta \sqrt{\lambda _{1}\left( \Delta _{f}\right) 
}\,r(A,B)}dt\,dz\,dy \\
&=&\frac{1}{2}\,\int_{A}\int_{B}G(y,z)\,e^{-f\left( y\right) }e^{-f\left(
z\right) }e^{\delta \sqrt{\lambda _{1}\left( \Delta _{f}\right) }%
\,r(A,B)}dz\,dy\,.
\end{eqnarray*}%
Combining with the previous inequality, we conclude

\begin{eqnarray}
&&\int_{A}\int_{B}G(y,z)e^{-f\left( y\right) }e^{-f\left( z\right) }\,dz\,dy
\label{b1} \\
&\leq &\frac{1}{\left( 1-\delta ^{2}\right) \lambda _{1}\left( \Delta
_{f}\right) }\sqrt{\mathrm{V}_{f}\left( A\right) }\sqrt{\mathrm{V}_{f}\left(
B\right) }e^{-\delta \,\sqrt{\lambda _{1}\left( \Delta _{f}\right) }%
\,r(A,B)}.  \notag
\end{eqnarray}%
Clearly, this proves the theorem if $r\left( A,B\right) \leq 1.$ When $%
r\left( A,B\right) >1,$ the theorem follows by setting

\begin{equation*}
\delta :=1-\frac{1}{r\left( A,B\right) }
\end{equation*}%
in (\ref{b1}).
\end{proof}

One may wish to compare Theorem \ref{GAB} with a result in \cite{LW} that

\begin{eqnarray}
&&\int_{B\left( x,R+1\right) \backslash B\left( x,R\right) }G^{2}\left(
x,y\right) e^{-f\left( y\right) }dy  \label{LW} \\
&\leq &Ce^{-2\sqrt{\lambda _{1}\left( \Delta _{f}\right) }R}\int_{B\left(
x,2\right) \backslash B\left( x,1\right) }G^{2}\left( x,y\right) e^{-f\left(
y\right) }dy.  \notag
\end{eqnarray}

Before we continue, let us first recall some results from \cite{MW2}
concerning smooth metric measure space $\left( M,g,e^{-f}dx\right) .$
Suppose that the associated Bakry-Emery Ricci curvature is bounded below by

\begin{equation}
\mathrm{Ric}_{f}\geq -\left( n-1\right) K  \label{Ric}
\end{equation}%
and weight $f$ satisfies

\begin{equation}
\sup_{y\in B\left( x,1\right) }\left\vert f\left( x\right) -f\left( y\right)
\right\vert \leq a.  \label{f}
\end{equation}%
Then Sobolev inequality of the following form holds.

\begin{equation}
\left( \int_{B\left( x,1\right) }\phi ^{\frac{2n}{n-2}}e^{-f}\right) ^{\frac{%
n-2}{n}}\leq \frac{c\left( n,K,a\right) }{\mathrm{V}_{f}\left( x,1\right) ^{%
\frac{2}{n}}}\left( \int_{B\left( x,1\right) }\left\vert \nabla \phi
\right\vert ^{2}e^{-f}+\int_{B\left( x,1\right) }\phi ^{2}e^{-f}\right)
\label{s}
\end{equation}%
for any $\phi $ with support in $B\left( x,1\right).$ Also, volume
comparison of the form

\begin{equation}
\frac{\mathrm{V}_{f}\left( x,r_{2}\right) }{\mathrm{V}_{f}\left(
x,r_{1}\right) }\leq c\left( a\right) \left( \frac{\int_{0}^{r_{2}}\sinh
\left( \sqrt{K}t\right) dt}{\int_{0}^{r_{1}}\sinh \left( \sqrt{K}t\right) dt}%
\right) ^{n-1+2a}  \label{v}
\end{equation}%
is valid for any $0<r_{1}<r_{2}<1.$ In view of (\ref{f}), both inequalities
hold with respect to the Riemannian volume as well. One also has the
gradient estimate of the form

\begin{equation}
\left\vert \nabla \ln u\right\vert \leq C\left( n,K,a\right) \text{ on }%
B\left( x,\frac{1}{20}\right)  \label{b2}
\end{equation}%
for any $u>0$ with $\Delta _{f}u=0$ in $B\left( x,\frac{1}{10}\right) .$

As a consequence of Theorem \ref{GAB}, we have the following.

\begin{proposition}
\label{Gball}Let $\left( M^n,\,g,\,e^{-f}dx\right) $ be a smooth metric
measure space satisfying (\ref{Ric}) and (\ref{f}). If $\lambda _{1}\left(
\Delta _{f}\right) >0,$ then the minimal positive Green's function $G(x,y)$
of $\Delta _{f}$ satisfies

\begin{equation*}
\int_{B\left( x,1\right) }G(x,y)\,e^{-f\left( y\right) }dy\leq C
\end{equation*}%
with $C$ depending on $n$, $K$, $a$ and $\lambda _{1}\left( \Delta
_{f}\right).$
\end{proposition}

\begin{proof}
Applying Theorem \ref{GAB} to $A=B=B\left( x,1\right) $ we get

\begin{equation*}
\int_{B\left( x,1\right) }\int_{B\left( x,1\right) }G(y,z)\,e^{-f\left(
y\right) }e^{-f\left( z\right) }dy\,dz\leq C\,\mathrm{V}_{f}\left(
x,1\right).
\end{equation*}%
Note that the function 
\begin{equation*}
u\left( z\right) :=\int_{B\left( x,1\right) }G(y,z)\,e^{-f\left( y\right) }dy
\end{equation*}%
satisfies $\Delta _{f}u=-1$ on $B\left( x,1\right) $ and 
\begin{equation*}
\int_{B\left( x,1\right) }u\left( z\right) e^{-f\left( z\right) }dz\leq C\,%
\mathrm{V}_{f}\left( x,1\right).
\end{equation*}%
Since $u>0,$ we have 
\begin{equation*}
\Delta _{f}\left( u+1\right) \geq -\left( u+1\right) \text{ on }B\left(
x,1\right).
\end{equation*}%
Using (\ref{s}) and (\ref{v}), the standard DeGiorgi-Nash-Moser iteration
implies that 
\begin{eqnarray*}
u\left( x\right) +1 &\leq &C\mathrm{V}_{f}^{-1}\left( x,1\right)
\int_{B\left( x,1\right) }\left( u\left( y\right) +1\right) e^{-f\left(
y\right) }dy \\
&\leq &C.
\end{eqnarray*}%
This proves the proposition.
\end{proof}

For the following proposition, we will work with the level sets of the
Green's function. Denote by

\begin{equation*}
L_{x}\left( \alpha ,\beta \right) :=\left\{ y\in M:\alpha <G\left(
x,y\right) <\beta \right\}.
\end{equation*}%
Also, we will use $c$ and $C$ to denote constants depending
on $n,$ $K,$ $a$ and $\lambda _{1}\left( \Delta _{f}\right).$ These constants
may change from line to line.

\begin{proposition}
\label{MS}Let $\left( M,\,g,\,e^{-f}dx\right) $ be a smooth metric measure
space satisfying (\ref{Ric}) and (\ref{f}). If $\lambda _{1}\left( \Delta
_{f}\right) >0,$ then

(i) for any $r>0,$%
\begin{equation*}
\sup_{y\in B\left( p,r\right) \backslash B\left( x,1\right) }G\left(
x,y\right) \leq e^{c\,r}\inf_{y\in B\left( p,r\right) \backslash
B\left( x,1\right) }G\left( x,y\right);
\end{equation*}

(ii) for $x\in M$ and $0<\alpha <\beta,$ 
\begin{equation*}
\int_{L_{x}\left( \alpha ,\beta \right) }G\left( x,y\right) e^{-f\left(
y\right) }dy\leq c\left( 1+\ln \frac{\beta }{\alpha }\right).
\end{equation*}
\end{proposition}

\begin{proof} For (i), we first show that 
\begin{equation}
G\left( x,y\right) \leq c\,G\left( x,z\right)  \label{b3}
\end{equation}%
for $y,z\in \partial B\left( x,1\right).$

For $y\in \partial B\left( x,1\right),$ since the function $G\left( x,\cdot
\right) $ is $f$-harmonic on $B\left( y,1\right),$ by (\ref{b2}),

\begin{equation*}
G\left( x,y\right) \leq c\,G\left( x,z\right) 
\end{equation*}%
for $z\in B\left( y,\frac{1}{5}\right).$ Hence, it suffices to prove (\ref{b3}) 
for $y$ and $z$ satisfying $d\left( y,z\right) \geq \frac{1}{5}.$ 
Let $\gamma \left(t\right) $ and $\sigma \left(t\right) $ be minimizing geodesics from $x$ to $y$ 
and from $x$ to $z$ respectively, $t\in \left[ 0,1\right].$ We have that $d\left( y,\sigma \right) \geq \frac{1}{10}$ 
and $d\left( z,\gamma \right) \geq \frac{1}{10}.$ Indeed, suppose that there exists 
$t_{0}\in \left( 0,1\right) $ such that $d\left( y,\sigma \left( t_{0}\right)
\right) <\frac{1}{10}.$ Since $d\left( x,y\right) =d\left( x,z\right) =1$
and $d\left( y,z\right) \geq \frac{1}{5},$ the triangle inequality implies 
\begin{eqnarray*}
d\left( z,\sigma \left( t_{0}\right) \right)  &\geq &d\left( y,z\right)
-d\left( y,\sigma \left( t_{0}\right) \right)  \\
&>&\frac{1}{10}
\end{eqnarray*}%
and 
\begin{eqnarray*}
d\left( x,\sigma \left( t_{0}\right) \right)  &\geq &d\left( x,y\right)
-d\left( y,\sigma \left( t_{0}\right) \right)  \\
&>&\frac{9}{10}.
\end{eqnarray*}%
Adding up these two inequalities we get 
\begin{eqnarray*}
d\left( x,z\right)  &=&d\left( x,\sigma \left( t_{0}\right) \right) +d\left(
\sigma \left( t_{0}\right) ,z\right)  \\
&>&1.
\end{eqnarray*}%
This contradiction shows that $d\left( y,\sigma \right) \geq \frac{1}{10}$ as claimed. 
The proof of $d\left( z,\gamma \right) \geq \frac{1}{10}$
is similar.

Consequently, $G\left( y,\cdot \right) $ is $f$-harmonic on $B\left( \sigma
\left( t\right) ,\frac{1}{10}\right) $ for all $t\in \left[ 0,1\right].$ It
follows from (\ref{b2}) that 
\begin{equation}
G\left( y,x\right) \leq c\,G\left( y,z\right).  \label{b3.1}
\end{equation}%
Similarly, as $d\left( z,\gamma \right) \geq \frac{1}{10},$ 
$G\left( z,\cdot \right) $ is $f$-harmonic on $B\left( \gamma \left( t\right)
\right) $ for $t\in \left[ 0,1\right].$ By (\ref{b2}) we get 
\begin{equation}
G\left( z,y\right) \leq c\,G\left( z,x\right).  \label{b3.2}
\end{equation}%
Combining (\ref{b3.1}) with (\ref{b3.2}) we conclude that 
$G\left( x,y\right) \leq c\,G\left( x,z\right) $ as claimed in (\ref{b3}).

Now for given $r>0,$ suppose first that $r>\frac{1}{2}d\left( p,x\right).$ For
$y,z\in B\left( p,r\right) \backslash B\left( x,1\right),$ let $\tau \left( t\right) $ 
and $\eta (t)$ be minimizing normal geodesics from $x$ to $y$ and from $x$ to $z$ respectively. 
We denote $y_{1}:=\tau \left( 1\right) $ and $z_{1}:=\eta \left( 1\right).$ Since $%
y_{1},z_{1}\in \partial B\left( x,1\right),$ by (\ref{b3}) we have 
 
\begin{equation}
G\left( x,y_{1}\right) \leq c\,G\left( x,z_{1}\right).  \label{b3.3}
\end{equation}%
On the other hand, the function $G\left( x,\cdot \right) $ is $f$-harmonic
on $B\left( \gamma \left( t\right) ,\frac{1}{10}\right) $ for all $t\geq 1.$
By (\ref{b2}) we obtain 

\begin{equation*}
G\left( x,y\right) \leq e^{cd\left( x,y\right) }\,G\left( x,y_{1}\right).
\end{equation*}%
Similarly, we have 
\begin{equation*}
G\left( x,z_{1}\right) \leq e^{cd\left( x,z\right) }\,G\left( x,z\right).
\end{equation*}%
In view of (\ref{b3.3}) we conclude that 
\begin{equation}
G\left( x,y\right) \leq e^{c\left( d\left( x,y\right) +d\left( x,z\right)
\right) }\,G\left( x,z\right).  \label{b3.4}
\end{equation}%
Since $d\left( p,x\right) <2r$ and $y,z\in B\left( p,r\right),$ 
by the triangle inequality, $d\left( x,y\right) <3r$ and $d\left(
x,z\right) <3r.$ Hence, (\ref{b3.4}) implies that 
\begin{equation*}
G\left( x,y\right) \leq e^{cr}\,G\left( x,z\right) 
\end{equation*}%
for $y,z\in B\left( p,r\right) \backslash B\left( x,1\right).$ This
proves (i) in the case that $r>\frac{1}{2}d\left( p,x\right).$

Suppose now that $r\leq \frac{1}{2}d\left( p,x\right).$ We
may assume that 
\begin{equation*}
B\left( p,r\right) \backslash B\left( x,1\right) \neq \varnothing.
\end{equation*}%
For $q\in B\left( p,r\right) \backslash B\left( x,1\right),$ we have 
$d\left( q,x\right) \geq 1$ and 
\begin{eqnarray*}
d\left( p,q\right)  &<&r \\
&\leq &\frac{1}{2}d\left( p,x\right).
\end{eqnarray*}%
The triangle inequality implies 
\begin{eqnarray*}
d\left( p,x\right)  &\geq &d\left( q,x\right) -d\left( p,q\right)  \\
&\geq &1-\frac{1}{2}d\left( p,x\right).
\end{eqnarray*}%
This shows that 
\begin{equation}
d\left( p,x\right) \geq \frac{2}{3}.  \label{b3.5}
\end{equation}%
We claim that 
\begin{equation}
x\notin B\left( p,r+\frac{1}{10}\right).  \label{b3.6}
\end{equation}%
Otherwise, $d\left( p,x\right) <r+\frac{1}{10}.$ Since $d\left( p,x\right) \geq 2r,$ 
it follows that $r<\frac{1}{10}$ and $d\left( p,x\right) <\frac{1}{5}.$ This contradicts with
(\ref{b3.5}). So (\ref{b3.6}) holds and $G\left( x,\cdot \right) $ is $f$-harmonic on $%
B\left( p,r+\frac{1}{10}\right).$  It is then easy to see from (\ref{b2}) that 
\begin{equation*}
\sup_{y\in B\left( p,r\right) }G\left( x,y\right) \leq e^{cr}\,\inf_{y\in
B\left( p,r\right) }G\left( x,y\right).
\end{equation*}%
This proves (i) in the remaining case that $r\leq \frac{1}{2}d\left(
p,x\right).$

To prove part (ii), let $\phi =\chi \psi $ be a cut-off function with
compact support on $M,$ where 
\begin{equation*}
\chi \left( y\right) :=\left\{ 
\begin{array}{c}
\ln \left( e\beta \right) -\ln G\,\left( x,y\right) \\ 
1 \\ 
\ln G\left( x,y\right) -\ln \left( e^{-1}\alpha \right) \\ 
0%
\end{array}%
\right. 
\begin{array}{l}
\text{on }L_{x}\left( \beta ,e\beta \right) \\ 
\text{on }L_{x}\left( \alpha ,\beta \right) \\ 
\text{on }L_{x}\left( e^{-1}\alpha ,\alpha \right) \\ 
\text{otherwise}%
\end{array}%
\end{equation*}%
and 
\begin{equation*}
\psi \left( y\right) =\left\{ 
\begin{array}{c}
1 \\ 
R+1-d\left( x,y\right) \\ 
0%
\end{array}%
\right. 
\begin{array}{l}
\text{on }B\left( x,R\right) \\ 
\text{on }B\left( x,R+1\right) \backslash B\left( x,R\right) \\ 
\text{on }M\backslash B\left( x,R+1\right)%
\end{array}%
\end{equation*}%
Obviously,

\begin{eqnarray}
&&\lambda _{1}\left( \Delta _{f}\right) \int_{M}\phi ^{2}\left( y\right)
G\left( x,y\right) e^{-f\left( y\right) }dy  \label{b4} \\
&\leq &\int_{M}\left\vert \nabla \left( \phi G^{\frac{1}{2}}\right)
\right\vert ^{2}\left( x,y\right) e^{-f\left( y\right) }dy  \notag \\
&\leq &\frac{1}{2}\int_{M}\phi ^{2}\left( y\right) \left\vert \nabla
G\right\vert ^{2}\left( x,y\right) G^{-1}\left( x,y\right) e^{-f\left(
y\right) }dy  \notag \\
&&+2\int_{M}G\left( x,y\right) \left\vert \nabla \phi \right\vert ^{2}\left(
y\right) e^{-f\left( y\right) }dy.  \notag
\end{eqnarray}%
To compute the integrals on the right side of (\ref{b4}) we use the co-area
formula 
\begin{eqnarray}
&&\int_{L_{x}\left( t_{0},t_{1}\right) }\left\vert \nabla G\right\vert
^{2}\left( x,y\right) G^{-1}\left( x,y\right) e^{-f\left( y\right) }dy
\label{b5} \\
&=&\int_{t_{0}}^{t_{1}}t^{-1}\left( \int_{l_{x}\left( t\right) }\left\vert
\nabla G\right\vert \left( x,y\right) e^{-f\left( y\right) }dy\right) dt, 
\notag
\end{eqnarray}%
where $0<t_{0}<t_{1}$ and 
\begin{equation*}
l_{x}\left( t\right) :=\left\{ y\in M:G\left( x,y\right) =t\right\}.
\end{equation*}

Note that although $L_{x}\left( t_{0},t_{1}\right) $ may not be compact,
both sides of (\ref{b5}) are finite with the equality justified in \cite{LW1}. 
Furthermore, 
\begin{equation}
\int_{l_{x}\left( t\right) }\left\vert \nabla G\right\vert \left( x,y\right)
e^{-f\left( y\right) }dy=1  \label{b6}
\end{equation}%
for all $t>0$ and for all $x\in M.$ Indeed, using that $G\left( x,y\right) $
is $f$-harmonic for $y\neq x,$ we have that 
\begin{equation*}
\int_{l_{x}\left( t\right) }\left\vert \nabla G\right\vert \left( x,y\right)
e^{-f\left( y\right) }dy=\int_{\partial B\left( x,\varepsilon \right) }\frac{%
\partial G}{\partial r}\left( x,y\right) e^{-f\left( y\right) }dy
\end{equation*}%
for sufficiently small $\varepsilon >0.$ Letting $\varepsilon \rightarrow 0$
and using the asymptotics of $G\left( x,y\right) $ as $y\rightarrow x,$ we
get (\ref{b6}). Hence, (\ref{b5}) becomes 
\begin{equation}
\int_{L_{x}\left( t_{0},t_{1}\right) }\left\vert \nabla G\right\vert
^{2}\left( x,y\right) G^{-1}\left( x,y\right) e^{-f\left( y\right) }dy=\ln 
\frac{t_{1}}{t_{0}}.  \label{b6'}
\end{equation}%
This implies that 
\begin{eqnarray}
&&\int_{M}\phi ^{2}\left( y\right) \left\vert \nabla G\right\vert ^{2}\left(
x,y\right) G^{-1}\left( x,y\right) e^{-f\left( y\right) }dy  \label{b7} \\
&\leq &\int_{L_{x}\left( e^{-1}\alpha ,e\beta \right) }\left\vert \nabla
G\right\vert ^{2}\left( x,y\right) G^{-1}\left( x,y\right) e^{-f\left(
y\right) }dy  \notag \\
&\leq &2+\ln \frac{\beta }{\alpha }.  \notag
\end{eqnarray}%
To estimate the second term of the right hand side of (\ref{b4}), note that

\begin{eqnarray}
&&\int_{M}G\left( x,y\right) \left\vert \nabla \phi \right\vert ^{2}\left(
y\right) e^{-f\left( y\right) }dy  \label{b8} \\
&\leq &2\int_{M}G\left( x,y\right) \left\vert \nabla \psi \right\vert
^{2}\left( y\right) \chi ^{2}\left( y\right) e^{-f\left( y\right) }dy  \notag
\\
&&+2\int_{M}G\left( x,y\right) \left\vert \nabla \chi \right\vert ^{2}\left(
y\right) \psi ^{2}\left( y\right) e^{-f\left( y\right) }dy.  \notag
\end{eqnarray}%
Since $G>e^{-1}\alpha $ on the support of $\chi,$ it follows that 
\begin{eqnarray}
&&\int_{M}G\left( x,y\right) \left\vert \nabla \psi \right\vert ^{2}\left(
y\right) \chi ^{2}\left( y\right) e^{-f\left( y\right) }dy  \label{b9} \\
&\leq &\frac{1}{\alpha }\int_{B\left( x,R+1\right) \backslash B\left(
x,R\right) }G^{2}\left( x,y\right) e^{-f\left( y\right) }dy  \notag \\
&\leq &\frac{C}{\alpha }e^{-2\sqrt{\lambda _{1}\left( \Delta _{f}\right) }%
R}\int_{B\left( x,2\right) \backslash B\left( x,1\right) }G^{2}\left(
x,y\right) e^{-f\left( y\right) }dy,  \notag
\end{eqnarray}%
where in the last line we have used (\ref{LW}). Furthermore, (\ref{b6'})
yields 
\begin{eqnarray*}
&&\int_{M}G\left( x,y\right) \left\vert \nabla \chi \right\vert ^{2}\left(
y\right) \psi ^{2}\left( y\right) e^{-f\left( y\right) }dy \\
&\leq &\int_{L_{x}\left( \beta ,e\beta \right) }\left\vert \nabla
G\right\vert ^{2}\left( x,y\right) G^{-1}\left( x,y\right) e^{-f\left(
y\right) }dy \\
&&+\int_{L_{x}\left( e^{-1}\alpha ,\alpha \right) }\left\vert \nabla
G\right\vert ^{2}\left( x,y\right) G^{-1}\left( x,y\right) e^{-f\left(
y\right) }dy \\
&\leq &c.
\end{eqnarray*}%
Combining this with (\ref{b9}) and (\ref{b8}) we obtain%
\begin{eqnarray*}
&&\int_{M}G\left( x,y\right) \left\vert \nabla \phi \right\vert ^{2}\left(
y\right) e^{-f\left( y\right) }dy \\
&\leq &\frac{C}{\alpha }e^{-2\sqrt{\lambda _{1}\left( \Delta _{f}\right) }%
R}\int_{B\left( x,2\right) \backslash B\left( x,1\right) }G^{2}\left(
x,y\right) e^{-f\left( y\right) }dy+C.
\end{eqnarray*}%
Together with (\ref{b7}) and (\ref{b4}), this implies that 
\begin{eqnarray*}
&&\lambda _{1}\left( \Delta _{f}\right) \int_{L_{x}\left( \alpha ,\beta
\right) \cap B\left( x,R\right) }G\left( x,y\right) e^{-f\left( y\right) }dy
\\
&\leq &\lambda _{1}\left( \Delta _{f}\right) \int_{M}\phi ^{2}\left(
y\right) G\left( x,y\right) e^{-f\left( y\right) }dy \\
&\leq &\frac{C}{\alpha }e^{-2\sqrt{\lambda _{1}\left( \Delta _{f}\right) }%
R}\int_{B\left( x,2\right) \backslash B\left( x,1\right) }G^{2}\left(
x,y\right) e^{-f\left( y\right) }dy \\
&&+\ln \frac{\beta }{\alpha }+C.
\end{eqnarray*}%
Letting $R\rightarrow \infty $ we conclude that%
\begin{equation*}
\int_{L_{x}\left( \alpha ,\beta \right) }G\left( x,y\right) e^{-f\left(
y\right) }dy\leq c\left( 1+\ln \frac{\beta }{\alpha }\right).
\end{equation*}%
So the proposition is proved.
\end{proof}

We now come to the following crucial estimate. The estimate is sharp as it
can be readily checked on the hyperbolic space $\mathbb{H}^{n}$ with trivial
weight $f.$

\begin{theorem}
\label{Gr}Let $\left( M,g,e^{-f}dx\right) $ be a smooth metric measure space
satisfying (\ref{Ric}) and (\ref{f}). If $\lambda _{1}\left( \Delta
_{f}\right) >0,$ then for fixed $p\in M$ and $r>0,$

\begin{equation*}
\int_{B\left( p,r\right) }G\left( x,y\right) e^{-f\left( y\right) }dy\leq
C\,\left( 1+r\right)
\end{equation*}%
for some constant $C$ depending on $n,$ $K,$ $a$ and $\lambda _{1}\left( \Delta
_{f}\right).$ 
\end{theorem}

\begin{proof}
We first prove that

\begin{equation}
\int_{B\left( p,r\right) \backslash B\left( x,1\right) }G\left( x,y\right)
e^{-f\left( y\right) }dy\leq C\left( 1+r\right) .  \label{b10}
\end{equation}%
Let

\begin{equation*}
\alpha :=\inf_{y\in B\left( p,r\right) \backslash B\left( x,1\right) }G\left(
x,y\right) \text{ and }\beta :=\sup_{y\in B\left( p,r\right) \backslash
B\left( x,1\right) }G\left( x,y\right) .
\end{equation*}%
It follows from part (ii) of Proposition \ref{MS} that 
\begin{eqnarray*}
\int_{B\left( p,r\right) \backslash B\left( x,1\right) }G\left( x,y\right)
e^{-f\left( y\right) }dy &\leq &\int_{L_{x}\left( \alpha ,\beta \right)
}G\left( x,y\right) e^{-f\left( y\right) }dy \\
&\leq &c\left( \ln \frac{\beta }{\alpha }+1\right) .
\end{eqnarray*}%
However, part (i) of Proposition \ref{MS} implies that 
\begin{equation*}
\beta \leq e^{c\,r}\alpha.
\end{equation*}%
Therefore, (\ref{b10}) follows. In view of Proposition \ref{Gball}, one concludes 
that
\begin{equation*}
\int_{B\left( p,r\right) }G\left( x,y\right) e^{-f\left( y\right) }dy\leq
C\,\left( 1+r\right)
\end{equation*}%
for any $r>0.$ This proves the theorem.
\end{proof}

Finally, we point out that 

\begin{equation}
\int_{\partial B\left( p,t\right) }G\left( x,y\right) e^{-f\left( y\right)
}dy\leq c  \label{b12}
\end{equation}%
for all $x\in M$ and $0<t\leq \frac{1}{2},$ where
$c$ is a constant depending on $n,$ $K$, $a,$ $\lambda_1(\Delta_f)$ and possibly the geometry
of $B(p,1).$ Indeed, if $x\in B\left(
p,1\right) $ then (\ref{b12}) is clearly true as $c$ is allowed to depend on
the geometry of $B\left( p,1\right).$ In the case of $d\left( p,x\right)
\geq 1,$ since $G\left( x,\cdot \right) $ is $f$-harmonic on $B\left(
p,1\right),$ by (\ref{b2}) we have

\begin{equation*}
\sup_{y\in B\left( p,\frac{1}{2}\right) }G\left( x,y\right) \leq c\inf_{y\in
B\left( p,\frac{1}{2}\right) }G\left( x,y\right).
\end{equation*}
Note by Theorem \ref{Gr},

\begin{equation*}
\int_{B\left( p,\frac{1}{2}\right) }G\left( x,y\right) e^{-f\left( y\right)
}dy\leq c.
\end{equation*}
It follows that $\inf_{y\in B\left(p,\frac{1}{2}\right) }G\left( x,y\right)
\leq c.$ Therefore, $\sup_{y\in B\left( p,\frac{1}{2}\right) }G\left(
x,y\right) \leq c$ as well. It is then easy to see that (\ref{b12}) is
indeed true.

Whether (\ref{b12}) is true for all $t>0$ remains an open question.

\section{Solving Poisson equation}

In this section, we solve the Poisson equation. We continue to denote by $%
\left( M,g,e^{-f}dx\right) $ an $n$-dimensional smooth metric measure space
satisfying the assumptions that

\begin{equation}
\lambda _{1}\left( \Delta _{f}\right) >0,  \label{s1}
\end{equation}%
the Bakry-Emery Ricci tensor is bounded below 
\begin{equation}
\mathrm{Ric}_{f}\geq -\left( n-1\right) K  \label{s2}
\end{equation}%
and the oscillation of $f$ on any unit ball $B\left( x,1\right) \subset M$
is uniformly bounded above by some fixed constant 
\begin{equation}
\sup_{y\in B\left( x,1\right) }\left\vert f\left( y\right) -f\left( x\right)
\right\vert \leq a.  \label{s3}
\end{equation}%
For simplicity, write $r\left( x\right) :=d\left( p,x\right),$ where $p\in M$
is a fixed point. We use $c,$ $C$ and $C_0$ to denote positive constants
depending on $n$, $K$, $a$, $\lambda _{1}\left( \Delta _{f}\right) $ and
possibly the geometry of $B\left( p,1\right).$

\begin{lemma}
\label{Vol}Let $\left( M,g,e^{-f}dx\right) $ be a smooth metric measure
space satisfying (\ref{s2}) and (\ref{s3}). Then the weighted volume
satisfies

\begin{equation*}
\mathrm{V}_{f}\left( p,R\right) \leq c\,R^{n+2a}\,e^{\left( \left(
n-1+2a\right) \sqrt{K}+a\right) R}
\end{equation*}%
for all $R\geq 1.$
\end{lemma}

\begin{proof}
Let us denote the volume form in geodesic coordinates centered at $p$ by 
\begin{equation*}
dV|_{\exp _{p}\left( r\xi \right) }=J\left( p,r,\xi \right) drd\xi
\end{equation*}%
for $r>0$ and $\xi \in S_{p}M,$ the unit tangent sphere at $p.$ Let $\gamma
\left( s\right) $ be a minimizing normal geodesic with $\gamma \left(
0\right) =p.$ Along $\gamma,$ according to the Laplace comparison theorem 
\cite{WW} we have 
\begin{eqnarray}
m_{f}\left( r\right) &\leq &\left( n-1\right) \sqrt{K}\coth \left( \sqrt{K}%
r\right)  \label{s5} \\
&&+\frac{2K}{\sinh ^{2}\left( \sqrt{Kr}\right) }\int_{0}^{r}\left( f\left(
s\right) -f\left( r\right) \right) \cosh \left( 2\sqrt{K}s\right) ds,  \notag
\end{eqnarray}%
where $m_{f}\left( r\right) :=\frac{d}{dr}\ln J_{f}\left( p,r,\xi \right) $
and $f\left( t\right) :=f\left( \gamma \left( t\right) \right) $. Using (\ref%
{s3}) we have that $\left\vert f\left( s\right) -f\left( r\right)
\right\vert \leq a\left( r-s+1\right).$ It follows that%
\begin{eqnarray*}
&&\int_{0}^{r}\left( f\left( s\right) -f\left( r\right) \right) \cosh \left(
2\sqrt{K}s\right) ds \\
&\leq &\frac{a}{2\sqrt{K}}\int_{0}^{r}\sinh \left( 2\sqrt{K}s\right) ds+%
\frac{a}{2\sqrt{K}}\sinh \left( 2\sqrt{K}r\right) \\
&=&\frac{a}{2K}\sinh ^{2}\left( \sqrt{K}r\right) +\frac{a}{2\sqrt{K}}\sinh
\left( 2\sqrt{K}r\right).
\end{eqnarray*}%
Therefore, we get from (\ref{s5}) that 
\begin{equation*}
m_{f}\left( r\right) \leq \left( \left( n-1+2a\right) \sqrt{K}\right) \coth
\left( \sqrt{K}r\right) +a.
\end{equation*}%
Thus, after integrating with respect to $r,$ 
\begin{equation*}
J_{f}\left( p,r,\xi \right) \leq r^{n-1+2a}e^{\left( \left( n-1+2a\right) 
\sqrt{K}+a\right) r}J_{f}\left( p,1,\xi \right).
\end{equation*}%
Integrating in $\xi \in S_{p}M$ then shows that 
\begin{equation*}
\mathrm{A}_{f}\left( p,r\right) \leq c\,r^{n-1+2a}\,e^{\left( \left(
n-1+2a\right) \sqrt{K}+a\right) r}.
\end{equation*}%
This proves the result.
\end{proof}

As $B\left( x,1\right) \subset B\left( p,r\left( x\right) +1\right),$ an
immediate consequence of Lemma \ref{Vol} is that 
\begin{equation}
\mathrm{V}_{f}\left( x,1\right) \leq c\,\left( 1+r\left( x\right) \right)
^{n+2a}\,e^{\left( \left( n-1+2a\right) \sqrt{K}+a\right) r\left( x\right) }
\label{Vx}
\end{equation}%
for all $x\in M.$

We are now ready to prove the main result of this paper. In the following, 
$\alpha_{0}$ is an arbitrary but fixed constant with 
\begin{equation}
0<\alpha_{0}<1.  \label{alpha}
\end{equation}

\begin{theorem}
\label{Poisson1}Let $\left( M,g,e^{-f}dx\right) $ be a smooth metric measure
space satisfying (\ref{s1}), (\ref{s2}) and (\ref{s3}). Let $\varphi $ be a
smooth function satisfying 
\begin{equation*}
\left\vert \varphi \right\vert \left( x\right) \leq \omega \left( r\left(
x\right) \right),
\end{equation*}%
where $\omega \left( t\right) $ is a non-increasing function such that $%
\int_{0}^{\infty }\omega \left( t\right) dt<\infty.$ Then the Poisson
equation $\Delta _{f}u=-\varphi $ admits a bounded solution $u$ on $M$ with

\begin{equation*}
\sup_{M}\left\vert u\right\vert \leq c\int_{0}^{\infty }\omega \left(
t\right) dt.
\end{equation*}%
Furthermore, there exists $C>0$ such that for all $\alpha \in [\alpha
_{0},1] $ and $x\in M,$

\begin{gather}
\left\vert u\right\vert \left( x\right) \leq C\left( \int_{2\alpha r\left(
x\right) }^{\infty }\omega \left( t\right) dt+r\left( x\right) \omega \left(
\alpha r\left( x\right) \right) \right)  \label{d} \\
+C\,\left( 1+r\left( x\right) \right)^{n+a}\,e^{-\sqrt{\lambda _{1}\left(
\Delta _{f}\right) }r\left( x\right) }\,\mathrm{V}_{f}^{-\frac{1}{2}}\left(
x,1\right)\, \int_{0}^{\alpha r\left( x\right) }\omega \left( t\right)
e^{bt}dt,  \notag
\end{gather}%
where 
\begin{equation*}
b:=\sqrt{\lambda _{1}\left( \Delta _{f}\right) }+\frac{1}{2}\left( \left(
\left( n-1\right) +2a\right) \sqrt{K}+a\right).
\end{equation*}
\end{theorem}

\begin{proof}
We first prove that 
\begin{equation}
\int_{M}G\left( x,y\right) \left\vert \varphi \right\vert \left( y\right)
e^{-f\left( y\right) }dy\leq c\int_{0}^{\infty }\omega \left( t\right) dt
\label{s6}
\end{equation}%
for all $x\in M.$ Note that

\begin{gather*}
\int_{B\left( p,1\right) }G\left( x,y\right) \left\vert \varphi \right\vert
\left( y\right) e^{-f\left( y\right) }dy=\int_{B\left( p,\frac{1}{2}\right)
}G\left( x,y\right) \left\vert \varphi \right\vert \left( y\right)
e^{-f\left( y\right) }dy \\
+\int_{B\left( p,1\right) \backslash B\left( p,\frac{1}{2}\right) }G\left(
x,y\right) \left\vert \varphi \right\vert \left( y\right) e^{-f\left(
y\right) }dy.
\end{gather*}%
By (\ref{b12}) and the co-area formula,

\begin{eqnarray*}
&&\int_{B\left( p,\frac{1}{2}\right) }G\left( x,y\right) \left\vert \varphi
\right\vert \left( y\right) e^{-f\left( y\right) }dy \\
&\leq &\int_{0}^{\frac{1}{2}}\left( \int_{\partial B\left( p,t\right)
}G\left( x,y\right) e^{-f\left( y\right) }dy\right) \sup_{\partial B\left(
p,t\right) }\left\vert \varphi \right\vert dt \\
&\leq &c\,\int_{0}^{\frac{1}{2}}\omega \left( t\right) dt,
\end{eqnarray*}%
where we have used that $\sup_{\partial B\left( p,t\right) }\left\vert
\varphi \right\vert \leq \omega \left( t\right) .$ By Theorem \ref{Gr} we get

\begin{eqnarray*}
\int_{B\left( p,1\right) \backslash B\left( p,\frac{1}{2}\right) }G\left(
x,y\right) \left\vert \varphi \right\vert \left( y\right) e^{-f\left(
y\right) }dy &\leq &c\sup_{B\left( p,1\right) \backslash B\left( p,\frac{1}{2%
}\right) }\left\vert \varphi \right\vert \\
&\leq &c\,\omega \left( \frac{1}{2}\right) \\
&\leq &c\,\int_{0}^{\frac{1}{2}}\omega \left( t\right) dt
\end{eqnarray*}%
as $\omega $ is non-increasing. In conclusion,

\begin{equation}
\int_{B\left( p,1\right) }G\left( x,y\right) \left\vert \varphi \right\vert
\left( y\right) e^{-f\left( y\right) }dy\leq c\int_{0}^{\frac{1}{2}}\omega
\left( t\right) dt.  \label{s6'}
\end{equation}%
Therefore,

\begin{eqnarray}
&&\int_{M}G\left( x,y\right) \left\vert \varphi \right\vert \left( y\right)
e^{-f\left( y\right) }dy  \label{s7} \\
&=&\sum_{j=0}^{\infty }\int_{B\left( p,2^{j+1}\right) \backslash B\left(
p,2^{j}\right) }G\left( x,y\right) \left\vert \varphi \right\vert \left(
y\right) e^{-f\left( y\right) }dy  \notag \\
&&+\int_{B\left( p,1\right) }G\left( x,y\right) \left\vert \varphi
\right\vert \left( y\right) e^{-f\left( y\right) }dy  \notag \\
&\leq &\sum_{j=0}^{\infty }\left( \int_{B\left( p,2^{j+1}\right) \backslash
B\left( p,2^{j}\right) }G\left( x,y\right) e^{-f\left( y\right) }dy\right)
\sup_{B\left( p,2^{j+1}\right) \backslash B\left( p,2^{j}\right) }\left\vert
\varphi \right\vert  \notag \\
&&+c\int_{0}^{\frac{1}{2}}\omega \left( t\right) dt.  \notag
\end{eqnarray}%
The hypothesis on $\varphi $ implies

\begin{equation*}
\sup_{B\left( p,2^{j+1}\right) \backslash B\left( p,2^{j}\right) }\left\vert
\varphi \right\vert \leq \omega \left( 2^{j}\right)
\end{equation*}%
and Theorem \ref{Gr} says that 
\begin{equation*}
\int_{B\left( p,2^{j+1}\right) \backslash B\left( p,2^{j}\right) }G\left(
x,y\right) e^{-f\left( y\right) }dy\leq c\,2^{j-1}.
\end{equation*}%
Using these estimates in (\ref{s7}) we obtain 
\begin{eqnarray*}
\int_{M}G\left( x,y\right) \left\vert \varphi \right\vert \left( y\right)
e^{-f\left( y\right) }dy &\leq &c\int_{0}^{\frac{1}{2}}\omega \left(
t\right) dt+c\,\sum_{j=0}^{\infty }2^{j-1}\,\omega \left( 2^{j}\right) \\
&\leq &c\int_{0}^{\infty }\omega \left( t\right) dt+c\sum_{j=0}^{\infty
}\int_{2^{j-1}}^{2^{j}}\omega \left( t\right) dt \\
&\leq &c\int_{0}^{\infty }\omega \left( t\right) dt.
\end{eqnarray*}%
This proves (\ref{s6}). As $\int_{0}^{\infty }\omega \left( t\right)
dt<\infty ,$ it follows that the function 
\begin{equation*}
u\left( x\right) :=\int_{M}G\left( x,y\right) \varphi \left( y\right)
e^{-f\left( y\right) }dy
\end{equation*}%
is well defined, bounded on $M,$ and verifies $\Delta _{f}u=-\varphi .$
Furthermore, we have the estimate 
\begin{equation}
\sup_{M}\left\vert u\right\vert \leq c\int_{0}^{\infty }\omega \left(
t\right) dt.  \label{s7'}
\end{equation}%
This proves the first part of the theorem.

We now prove the decay estimate (\ref{d}). For $x\in M,$ denote 
\begin{equation*}
R:=r\left( x\right) =d\left( p,x\right).
\end{equation*}%
Given $\alpha \in [\alpha _{0},1],$ let us first assume that $\alpha \,R<4.$
It follows that $r\left( x\right) \leq C.$ By (\ref{Vx}), it is obvious that

\begin{gather*}
\int_{0}^{\infty }\omega \left( t\right) dt\leq C\left( \int_{2\alpha
r\left( x\right) }^{\infty }\omega \left( t\right) dt+r\left( x\right)
\omega \left( \alpha r\left( x\right) \right) \right) \\
+C\left( 1+r\left( x\right) \right) ^{n+a}e^{-\sqrt{\lambda _{1}\left(
\Delta _{f}\right) }r\left( x\right) }\mathrm{V}_{f}^{-\frac{1}{2}}\left(
x,1\right) \int_{0}^{\alpha r\left( x\right) }\omega \left( t\right)
e^{bt}dt.
\end{gather*}%
In view of (\ref{s7'}), this proves (\ref{d}) when $\alpha \,R<4.$

From now on we assume that $\alpha \,R\geq 4.$ Note that

\begin{eqnarray*}
&&\int_{M\backslash B\left( p,\alpha R\right) }G\left( x,y\right) \left\vert
\varphi \right\vert \left( y\right) e^{-f\left( y\right) }dy \\
&=&\sum_{j=0}^{\infty }\int_{B\left( p,2^{j+1}\alpha R\right) \backslash
B\left( p,2^{j}\alpha R\right) }G\left( x,y\right) \left\vert \varphi
\right\vert \left( y\right) e^{-f\left( y\right) }dy \\
&\leq &\sum_{j=0}^{\infty }\left( \int_{B\left( p,2^{j+1}\alpha R\right)
\backslash B\left( p,2^{j}\alpha R\right) }G\left( x,y\right) e^{-f\left(
y\right) }dy\right) \sup_{B\left( p,2^{j+1}\alpha R\right) \backslash
B\left( p,2^{j}\alpha R\right) }\left\vert \varphi \right\vert .
\end{eqnarray*}%
Using the decay hypothesis on $\varphi $ we get 
\begin{equation*}
\sup_{B\left( p,2^{j+1}\alpha R\right) \backslash B\left( p,2^{j}\alpha
R\right) }\left\vert \varphi \right\vert \leq \omega \left( 2^{j}\alpha
R\right) .
\end{equation*}%
We also infer from Theorem \ref{Gr} that

\begin{equation*}
\int_{B\left( p,2^{j+1}\alpha R\right) \backslash B\left( p,2^{j}\alpha
R\right) }G\left( x,y\right) e^{-f\left( y\right) }dy\leq c\,2^{j-1}\,\alpha
\,R.
\end{equation*}%
In conclusion,

\begin{equation}
\int_{M\backslash B\left( p,\alpha R\right) }G\left( x,y\right) \left\vert
\varphi \right\vert \left( y\right) e^{-f\left( y\right) }dy\leq
c\,\sum_{j=0}^{\infty }\left( 2^{j-1}\,\alpha \,R\right) \,\omega \left(
2^{j}\,\alpha R\right) .  \label{s8}
\end{equation}%
However,

\begin{eqnarray*}
\sum_{j=0}^{\infty }2^{j-1}\,\alpha R\,\omega \left( 2^{j}\alpha \,R\right)
&\leq &c\,R\,\omega \left( \alpha \,R\right) +\sum_{j=2}^{\infty }\left(
2^{j-1}\,\alpha \,R\right) \omega \left( 2^{j}\,\alpha \,R\right) \\
&\leq &c\,R\,\omega \left( \alpha \,R\right) +\sum_{j=2}^{\infty
}\int_{2^{j-1}\,\alpha \,R}^{2^{j}\,\alpha \,R}\omega \left( t\right) dt \\
&=&c\,R\,\omega \left( \alpha \,R\right) +\int_{2\alpha \,R}^{\infty }\omega
\left( t\right) dt.
\end{eqnarray*}%
It follows that

\begin{equation}
\int_{M\backslash B\left( p,\alpha \,R\right) }G\left( x,y\right) \left\vert
\varphi \right\vert \left( y\right) e^{-f\left( y\right) }dy\leq
c\,R\,\omega \left( \alpha \,R\right) +c\int_{2\alpha \,R}^{\infty }\omega
\left( t\right) dt.  \label{s9}
\end{equation}

We now proceed to obtain a similar estimate on $B\left( p,\alpha R\right).$
By Theorem \ref{GAB},

\begin{eqnarray*}
&&\int_{B\left( x,1\right) }\int_{B\left( p,j+1\right) \backslash B\left(
p,j\right) }G(z,y)\,e^{-f\left( z\right) }e^{-f\left( y\right) }dy\,dz \\
&\leq &c\,R\,\sqrt{\mathrm{V}_{f}\left( x,1\right) }\,\sqrt{\mathrm{V}%
_{f}\left( p,j+1\right) }\,e^{-\sqrt{\lambda _{1}\left( \Delta _{f}\right) }%
\left( R-j\right) }
\end{eqnarray*}%
for any $j\in \left\{ 0,1,...,\left[ R\right] -3\right\},$ where $\left[ R%
\right] $ denotes the greatest integer less than or equal to $R.$ Using
Lemma \ref{Vol} that 
\begin{equation*}
\sqrt{\mathrm{V}_{f}\left( p,j+1\right) }\leq c\,R^{\frac{n}{2}+a}\,e^{\frac{%
1}{2}\left( \left( n-1+2a\right) \sqrt{K}+a\right) j},
\end{equation*}%
we conclude that

\begin{eqnarray}
&&\int_{B\left( x,1\right) }\int_{B\left( p,j+1\right) \backslash B\left(
p,j\right) }G(z,y)\,e^{-f\left( z\right) }e^{-f\left( y\right) }dy\,dz
\label{s10} \\
&\leq &c\,R^{n+a}\,e^{-\sqrt{\lambda _{1}\left( \Delta _{f}\right) }R}\,%
\sqrt{\mathrm{V}_{f}\left( x,1\right) }\,e^{bj},  \notag
\end{eqnarray}%
where 
\begin{equation*}
b=\sqrt{\lambda _{1}\left( \Delta _{f}\right) }+\frac{1}{2}\left( \left(
\left( n-1\right) +2a\right) \sqrt{K}+a\right).
\end{equation*}%
Note that for any $j\leq \left[ R\right] -3,$ 
\begin{equation*}
B\left( x,2\right) \cap \left( B\left( p,j+1\right) \backslash B\left(
p,j\right) \right) =\varnothing.
\end{equation*}%
Hence the function

\begin{equation*}
H\left( z\right) :=\int_{B\left( p,j+1\right) \backslash B\left( p,j\right)
}G(z,y)\,e^{-f\left( y\right) }dy
\end{equation*}%
is $f$-harmonic on $B\left( x,2\right).$ Applying (\ref{b2}) we get that 
\begin{equation*}
H\left( x\right) \leq c\,\mathrm{V}_{f}^{-1}\left( x,1\right) \int_{B\left(
x,1\right) }H\left( z\right) e^{-f\left( z\right) }dz.
\end{equation*}%
Together with (\ref{s10}), this gives

\begin{equation}
\int_{B\left( p,j+1\right) \backslash B\left( p,j\right)
}G(x,y)\,e^{-f\left( y\right) }dy\leq c\,R^{n+a}\, e^{-\sqrt{%
\lambda_{1}\left( \Delta _{f}\right) }R}\mathrm{V}_{f}^{-\frac{1}{2}}\left(
x,1\right) e^{bj}  \label{s11}
\end{equation}%
for $0\leq j\leq \left[ R\right] -3.$ We claim that (\ref{s11}) holds for $%
\left[ R\right] -3\leq j\leq \left[ R\right] $ as well. Indeed, in this case
(\ref{s11}) is equivalent to

\begin{equation*}
\int_{B\left( p,j+1\right) \backslash B\left( p,j\right)
}G(x,y)\,e^{-f\left( y\right) }dy\leq cR^{n+a}\mathrm{V}_{f}^{-\frac{1}{2}%
}\left( x,1\right) e^{\frac{1}{2}\left( \left( n-1+2a\right) \sqrt{K}%
+a\right) R}.
\end{equation*}%
This follows from Theorem \ref{Gr} that 
\begin{equation*}
\int_{B\left( p,j+1\right) \backslash B\left( p,j\right)
}G(x,y)\,e^{-f\left( y\right) }dy\leq c\,R
\end{equation*}%
for $\left[ R\right] -3\leq j\leq \left[ R\right] $ together with (\ref{Vx}%
). In conclusion, (\ref{s11}) holds true for all $j\in \left\{ 0,1,\cdots,%
\left[ R\right] \right\}.$

Following a similar argument as in (\ref{b12}), we have

\begin{equation}
\int_{\partial B\left( p,t\right) }G\left( x,y\right) e^{-f\left( y\right)
}dy\leq cR^{n+a}e^{-\sqrt{{\lambda}_{1}\left( \Delta _{f}\right) }R}\mathrm{V%
}_{f}^{-\frac{1}{2}}\left( x,1\right)  \label{s12}
\end{equation}%
for all $t\in \left[ 0,1\right].$ Indeed, (\ref{s11}) implies that 
\begin{eqnarray*}
\inf_{B\left( p,1\right) }G\left( x,y\right) &\leq &c\int_{B\left(
p,1\right) }G(x,y)\,e^{-f\left( y\right) }dy \\
&\leq &c\,R^{n+a}\,e^{-\sqrt{{\lambda}_{1}\left( \Delta _{f}\right) }R}%
\mathrm{V}_{f}^{-\frac{1}{2}}\left( x,1\right).
\end{eqnarray*}%
As the function $G\left( x,\cdot \right) $ is $f$-harmonic on $B\left(
p,2\right),$ one sees that 
\begin{equation*}
\sup_{B\left( p,1\right) }G\left( x,y\right) \leq c\,R^{n+a}\,e^{-\sqrt{{%
\lambda}_{1}\left( \Delta _{f}\right) }R}\mathrm{V}_{f}^{-\frac{1}{2}}\left(
x,1\right).
\end{equation*}%
This immediately implies (\ref{s12}).

We now write 
\begin{eqnarray}
&&\int_{B\left( p,\alpha R\right) }G(x,y)\,\left\vert \varphi \right\vert
\left( y\right) e^{-f\left( y\right) }dy  \label{s13} \\
&\leq &\sum_{j=1}^{\left[ \alpha R\right] }\int_{B\left( p,j+1\right)
\backslash B\left( p,j\right) }G(x,y)\,\left\vert \varphi \right\vert \left(
y\right) e^{-f\left( y\right) }dy  \notag \\
&&+\int_{B\left( p,1\right) }G(x,y)\,\left\vert \varphi \right\vert \left(
y\right) e^{-f\left( y\right) }dy.  \notag
\end{eqnarray}%
Using (\ref{s11}) we get 
\begin{eqnarray*}
&&\sum_{j=1}^{\left[ \alpha R\right] }\int_{B\left( p,j+1\right) \backslash
B\left( p,j\right) }G(x,y)\,\left\vert \varphi \right\vert \left( y\right)
e^{-f\left( y\right) }dy \\
&\leq &c\,\sum_{j=1}^{\left[ \alpha R\right] }\left( R^{n+a}\,e^{-\sqrt{{%
\lambda }_{1}\left( \Delta _{f}\right) }R}\mathrm{V}_{f}^{-\frac{1}{2}%
}\left( x,1\right) e^{bj}\right) \sup_{B\left( p,j+1\right) \backslash
B\left( p,j\right) }\left\vert \varphi \right\vert \\
&\leq &c\,R^{n+a}\,e^{-\sqrt{{\lambda }_{1}\left( \Delta _{f}\right) }R}%
\mathrm{V}_{f}^{-\frac{1}{2}}\left( x,1\right) \sum_{j=1}^{\left[ \alpha R%
\right] }\omega \left( j\right) e^{bj} \\
&\leq &c\,R^{n+a}\,e^{-\sqrt{{\lambda }_{1}\left( \Delta _{f}\right) }R}%
\mathrm{V}_{f}^{-\frac{1}{2}}\left( x,1\right) \int_{0}^{\alpha R}\omega
\left( t\right) e^{bt}dt,
\end{eqnarray*}%
where we have used that $\sup_{\partial B\left( p,t\right) }\left\vert
\varphi \right\vert \leq \omega \left( t\right) $ and that $\omega $ is
non-increasing in $t.$ By (\ref{s12}),

\begin{eqnarray*}
&&\int_{B\left( p,1\right) }G(x,y)\,\left\vert \varphi \right\vert \left(
y\right) e^{-f\left( y\right) }dy \\
&\leq &\int_{0}^{1}\left( \int_{\partial B\left( p,t\right) }G\left(
x,y\right) e^{-f\left( y\right) }dy\right) \sup_{\partial B\left( p,t\right)
}\left\vert \varphi \right\vert dt \\
&\leq &c\,R^{n+a}\,e^{-\sqrt{{\lambda }_{1}\left( \Delta _{f}\right) }R}%
\mathrm{V}_{f}^{-\frac{1}{2}}\left( x,1\right) \int_{0}^{1}\omega \left(
t\right) dt.
\end{eqnarray*}%
Plugging these estimates in (\ref{s13}) and using that $\alpha R\geq 4,$ we
conclude

\begin{eqnarray}
&&\int_{B\left( p,\alpha R\right) }G(x,y)\,\left\vert \varphi \right\vert
\left( y\right) e^{-f\left( y\right) }dy  \label{s14} \\
&\leq &c\,R^{n+a}\,e^{-\sqrt{{\lambda }_{1}\left( \Delta _{f}\right) }R}%
\mathrm{V}_{f}^{-\frac{1}{2}}\left( x,1\right) \int_{0}^{\alpha R}\omega
\left( t\right) e^{bt}dt.  \notag
\end{eqnarray}%
Finally, combining (\ref{s14}) and (\ref{s9}) we arrive at

\begin{eqnarray*}
&&\int_{M}G(x,y)\,\left\vert \varphi \right\vert \left( y\right) e^{-f\left(
y\right) }dy\leq c\int_{2\alpha R}^{\infty }\omega \left( t\right)
dt+c\,R\,\omega \left( \alpha R\right) \\
&&+c\,R^{n+a}\,e^{-\sqrt{{\lambda }_{1}\left( \Delta _{f}\right) }R}\mathrm{V%
}_{f}^{-\frac{1}{2}}\left( x,1\right) \int_{0}^{\alpha R}\omega \left(
t\right) e^{bt}dt.
\end{eqnarray*}%
This proves the theorem.
\end{proof}

In the case that the function $\varphi $ decays as 
\begin{equation*}
\left\vert \varphi \right\vert \left( x\right) \leq c\left( 1+r\left(
x\right) \right) ^{-k}
\end{equation*}%
for some $k>1$ and the weighted volume of unit balls is uniformly bounded
from below 
\begin{equation*}
\mathrm{V}_{f}\left( x,1\right) \geq c>0
\end{equation*}%
for all $x\in M,$ Theorem \ref{Poisson1} implies that the solution $u$
satisfies

\begin{gather*}
\left\vert u\right\vert \left( x\right) \leq C\,\left( \int_{2\alpha r\left(
x\right) }^{\infty }\omega \left( t\right) dt+r\left( x\right) \omega \left(
\alpha r\left( x\right) \right) \right) \\
+C\,\left( 1+r\left( x\right) \right) ^{n+a}\,e^{-\sqrt{\lambda _{1}\left(
\Delta _{f}\right) }r\left( x\right) }\int_{0}^{\alpha r\left( x\right)
}\omega \left( t\right) e^{bt}dt,
\end{gather*}%
where $\omega \left( t\right) =c\left( 1+t\right) ^{-k}.$ Taking $\alpha
=\alpha _{0}:=\frac{1}{2}b^{-1}\sqrt{\lambda _{1}\left( \Delta _{f}\right) }$
and estimating 
\begin{eqnarray*}
\int_{0}^{\alpha r\left( x\right) }\omega \left( t\right) e^{bt}dt &\leq
&ce^{\alpha br\left( x\right) } \\
&\leq &ce^{\frac{1}{2}\sqrt{\lambda _{1}\left( \Delta _{f}\right) }r\left(
x\right) },
\end{eqnarray*}%
one concludes that 
\begin{equation*}
\left\vert u\right\vert \left( x\right) \leq C\left( k\right) \left(
1+r\left( x\right) \right) ^{-k+1}
\end{equation*}%
as claimed by Theorem \ref{intro}.

As an application of Theorem \ref{Poisson1} we prove the following.

\begin{theorem}
\label{Nonl}Let $\left( M,g,e^{-f}dx\right) $ be a smooth metric measure
space satisfying (\ref{s1}), (\ref{s2}) and (\ref{s3}). Assume that the
weighted volume has lower bound $\mathrm{V}_{f}\left( x,1\right) \geq c>0$
for all $x\in M$. Suppose $\psi \geq 0$ satisfies%
\begin{equation*}
\Delta _{f}\psi \geq -c\psi ^{q}
\end{equation*}%
for some $q>1$, and 
\begin{equation*}
\lim_{x\rightarrow \infty }\psi \left( x\right) r^{\frac{1}{q-1}}\left(
x\right) =0.
\end{equation*}%
Then there exists $\delta >0$ and $C>0$ such that 
\begin{equation*}
\psi \left( x\right) \leq C\,e^{-r^{\delta }\left( x\right) }.
\end{equation*}
\end{theorem}

\begin{proof}
We first prove that for $\phi \geq 0$ satisfying

\begin{equation*}
\Delta _{f}\phi \geq -c\phi ^{q}
\end{equation*}
and 
\begin{equation*}
\phi \left( x\right) \leq \sigma \left( r\left( x\right) \right) 
\end{equation*}%
for a non-increasing function $\sigma \left( t\right) $ with $\int_{0}^{\infty
}\sigma ^{q}\left( t\right) dt<\infty,$ there exists $C_{0}>0$ and $%
\beta >0$ such that 
\begin{equation}
\phi \left( x\right) \leq C_{0}\left( \int_{\beta r\left( x\right) }^{\infty
}\sigma ^{q}\left( t\right) dt+e^{-\beta r\left( x\right) }\sigma ^{q}\left(
0\right) \right).  \label{n2}
\end{equation}

Indeed, by Theorem \ref{Poisson1} the equation $\Delta _{f}u=-c\phi ^{q}$
has a solution $u\geq 0$ satisfying 
\begin{eqnarray*}
u\left( x\right) &\leq &C\int_{\frac{\alpha }{2}r\left( x\right) }^{\infty
}\sigma ^{q}\left( t\right) dt \\
&&+C\left( 1+r\left( x\right) \right) ^{n+a}e^{-\sqrt{\lambda _{1}\left(
\Delta _{f}\right) }r\left( x\right) }\int_{0}^{\alpha r\left( x\right)
}\sigma ^{q}\left( t\right) e^{bt}dt
\end{eqnarray*}%
for all $x\in M.$ We let 
\begin{equation*}
\alpha :=\frac{1}{4}\frac{\sqrt{\lambda _{1}\left( \Delta _{f}\right) }}{b}
\end{equation*}%
and estimate 
\begin{eqnarray*}
\int_{0}^{\alpha r\left( x\right) }\sigma ^{q}\left( t\right) e^{bt}dt &\leq
&\sigma ^{q}\left( 0\right) \int_{0}^{\alpha r\left( x\right) }e^{bt}dt \\
&\leq &\frac{1}{b}\sigma ^{q}\left( 0\right) e^{\frac{1}{4}\sqrt{\lambda
_{1}\left( \Delta _{f}\right) }r\left( x\right) }.
\end{eqnarray*}%
It follows from above that 
\begin{equation}
u\left( x\right) \leq C\left( \int_{\frac{\alpha }{2}r\left( x\right)
}^{\infty }\sigma ^{q}\left( t\right) dt+e^{-\frac{1}{2}\sqrt{\lambda
_{1}\left( \Delta _{f}\right) }r\left( x\right) }\sigma ^{q}\left( 0\right)
\right).  \label{n3}
\end{equation}%
In particular, $u$ converges to zero at infinity. Since $\Delta _{f}\phi
\geq -c\phi ^{q}$ and $\Delta _{f}u=-c\phi ^{q},$ by the maximum principle
we get $\phi \leq u$ on $M.$ Therefore, (\ref{n2}) follows from (\ref{n3}) 
by setting

\begin{equation*}
\beta :=\min \left\{ \frac{\alpha }{2},\frac{1}{2}\sqrt{\lambda _{1}\left(
\Delta _{f}\right) }\right\}.
\end{equation*}

Now let 
\begin{equation*}
m_{0}=\left[ \left( q-1\right) ^{-2}\right] +2.
\end{equation*}%
Note that $q^{m}-m>0$ for $m\geq m_{0}.$
Fix $\varepsilon >0$ to be specified later. We prove by induction on $m\geq
m_{0}$ that 
\begin{equation}
\psi \left( x\right) \leq \varepsilon ^{q^{m}+m}\left( \beta ^{m}r\left(
x\right) +1\right) ^{-\frac{1}{q-1}}+B^{q^{m}-m}e^{-\beta ^{m}r\left(
x\right) },  \label{n6}
\end{equation}%
where $B$ is a large enough constant depending on $\varepsilon.$

First, note that (\ref{n6}) holds for $m=m_{0}$ by the assumption that 
\begin{equation*}
\lim_{x\rightarrow \infty }\psi \left( x\right) r^{\frac{1}{q-1}}\left(
x\right) =0
\end{equation*}%
and by adjusting the constant $B$ if necessary.

We now assume (\ref{n6}) holds for $m$ and prove 
\begin{equation}
\psi \left( x\right) \leq \varepsilon ^{q^{m+1}+\left( m+1\right) }\left(
\beta ^{m+1}r\left( x\right) +1\right) ^{-\frac{1}{q-1}}+B^{q^{m+1}-\left(
m+1\right) }e^{-\beta ^{m+1}r\left( x\right) }.  \label{n7}
\end{equation}

By the induction hypothesis we have $\psi \left( x\right) \leq \sigma \left(
r\left( x\right) \right),$ where 
\begin{equation*}
\sigma \left( t\right) :=\varepsilon ^{q^{m}+m}\left( \beta ^{m}t+1\right)
^{-\frac{1}{q-1}}+B^{q^{m}-m}e^{-\beta ^{m}t}
\end{equation*}%
is decreasing and $\int_{0}^{\infty }\sigma ^{q}\left( t\right) dt<\infty.$
Applying (\ref{n2}) we get that 
\begin{equation}
\psi \left( x\right) \leq C_{0}\left( \int_{\beta r\left( x\right) }^{\infty
}\sigma ^{q}\left( t\right) dt+e^{-\beta r\left( x\right) }\sigma ^{q}\left(
0\right) \right).  \label{n8}
\end{equation}%
Obviously, 

\begin{equation}
\sigma ^{q}\left( t\right) \leq c\varepsilon ^{q^{m+1}+qm}\left( \beta
^{m}t+1\right) ^{-\frac{q}{q-1}}+cB^{q^{m+1}-qm}e^{-q\cdot \beta ^{m}t}.
\label{n9}
\end{equation}%
It follows that%
\begin{eqnarray}
\int_{\beta r\left( x\right) }^{\infty }\sigma ^{q}\left( t\right) dt &\leq
&c\varepsilon ^{q^{m+1}+qm}\int_{\beta r\left( x\right) }^{\infty }\left(
\beta ^{m}t+1\right) ^{-\frac{q}{q-1}}dt  \label{n10} \\
&&+cB^{q^{m+1}-qm}\int_{\beta r\left( x\right) }^{\infty }e^{-q\cdot \beta
^{m}t}dt  \notag \\
&=&c\beta ^{-m}\varepsilon ^{q^{m+1}+qm}\left( \beta ^{m+1}r\left( x\right)
+1\right) ^{-\frac{1}{q-1}}  \notag \\
&&+c\beta ^{-m}B^{q^{m+1}-qm}e^{-\beta ^{m+1}r\left( x\right) }.  \notag
\end{eqnarray}%
Furthermore, we have by (\ref{n9}) that 
\begin{eqnarray}
e^{-\beta r\left( x\right) }\sigma ^{q}\left( 0\right) &\leq &c\left(
\varepsilon ^{q^{m+1}+qm}+B^{q^{m+1}-qm}\right) e^{-\beta r\left( x\right) }
\label{n11} \\
&\leq &c\beta ^{-m}B^{q^{m+1}-qm}e^{-\beta ^{m+1}r\left( x\right) }.  \notag
\end{eqnarray}%
Plugging (\ref{n10}) and (\ref{n11}) into (\ref{n8}) yields%
\begin{eqnarray}
\psi \left( x\right) &\leq &\left( cC_{0}\beta ^{-m}\varepsilon ^{qm-\left(
m+1\right) }\right) \varepsilon ^{q^{m+1}+\left( m+1\right) }\left( \beta
^{m+1}r\left( x\right) +1\right) ^{-\frac{1}{q-1}}  \label{n12} \\
&&+\left( cC_{0}\beta ^{-m}B^{-qm+\left( m+1\right) }\right)
B^{q^{m+1}-\left( m+1\right) }e^{-\beta ^{m+1}r\left( x\right) }.  \notag
\end{eqnarray}%
Since $m\geq m_{0}$, we have $qm-\left( m+1\right) \geq m+1.$ Now take 
$\varepsilon $ sufficiently small so that $\varepsilon \beta ^{-1}\leq 1$ and 
$cC_{0}\varepsilon \leq 1,$ and $B$ sufficiently large so
that $B^{-1}\beta ^{-1}\leq 1$ and $cC_{0}B^{-1}\leq 1.$ It follows by (\ref{n12}) that
\begin{equation*}
\psi \left( x\right) \leq \varepsilon ^{q^{m+1}+\left( m+1\right) }\left(
\beta ^{m+1}r\left( x\right) +1\right) ^{-\frac{1}{q-1}}+B^{q^{m+1}-\left(
m+1\right) }e^{-\beta ^{m+1}r\left( x\right) }.
\end{equation*}%
This proves (\ref{n7}). Hence,

\begin{equation}
\psi \left( x\right) \leq \varepsilon ^{q^{m}+m}\left( \beta ^{m}r\left(
x\right) +1\right) ^{-\frac{1}{q-1}}+B^{q^{m}-m}e^{-\beta ^{m}r\left(
x\right) }  \label{n13}
\end{equation}%
for all $m\geq m_{0}$ and $x\in M.$ 

For $x\in M$ fixed, we take 
\begin{equation*}
m:=\left[ \frac{\ln r\left( x\right) }{2\ln \left( q\beta ^{-1}\right) }%
\right],
\end{equation*}%
where $\left[ \cdot \right] $ denotes the greatest integer function. Here we
may assume that $r\left( x\right) $ is large enough so that $m\geq m_{0}.$
This implies that 
\begin{equation*}
B^{q^{m}-m}e^{-\beta ^{m}r\left( x\right) }\leq ce^{-c\sqrt{r}\left(
x\right) }
\end{equation*}%
and 
\begin{equation*}
\varepsilon ^{q^{m}+m}\left( \beta ^{m}r\left( x\right) +1\right) ^{-\frac{1%
}{q-1}}\leq ce^{-r^{\delta }\left( x\right) }
\end{equation*}%
for some $\delta >0$ depending on $\varepsilon.$ Hence, from (\ref{n13}) we conclude
that 

\begin{equation*}
\psi \left( x\right) \leq ce^{-r^{\delta }\left( x\right) }\text{ \ for all }%
x\in M.
\end{equation*}
\end{proof}

We conclude this section by showing the following simple proposition. Note
that the positivity of $\lambda _{1}(\Delta _{f})$ implies that the weighted
volume of $M$ must be infinite and the bottom spectrum $\lambda
_{1}^{M\backslash \Omega }\left( \Delta _{f}\right) $ of the weighted
Laplacian $\Delta _{f}$ on $M\backslash \Omega $ subject to the Dirichlet
boundary conditions on $\partial \Omega $ is positive as well. The
proposition says the converse is also true.

\begin{proposition}
Assume that $\lambda _{1}^{M\backslash \Omega }\left( \Delta _{f}\right) >0$
for compact domain $\Omega \subset M$ and the weighted volume of $M$ is
infinite. Then $M$ has positive spectrum $\lambda _{1}\left( \Delta
_{f}\right) >0.$
\end{proposition}

\begin{proof}
Pick $r_{0}>0$ so that $\Omega \subset B\left( p,r_{0}\right) .$ Since $%
\lambda _{1}^{M\backslash \Omega }\left( \Delta _{f}\right) >0,$ adapting
the results in \cite{LW, LW1} to $\Delta _{f}$ we know that $M$ is $f$%
-nonparabolic as it has infinite weighted volume. In particular, there
exists a positive nonconstant $f$-superharmonic, but not $f$-harmonic,
function $u$ on $M.$ We now claim that there also exists a positive,
strictly $f$-superharmonic function $w$ on $M.$ This can be seen as follows.
Consider the heat equation $\left( \Delta _{f}-\frac{\partial }{\partial t}%
\right) w=0$ with $w\left( 0\right) =u.$ Then $\left( \Delta _{f}-\frac{%
\partial }{\partial t}\right) \left( \Delta _{f}w\right) =0$ as well. By the
strong maximum principle, as $\Delta _{f}w\leq 0$ at $t=0,$ it follows that $%
\Delta _{f}w<0$ for $t>0.$ Hence, we obtain a positive, strictly $f$%
-superharmonic function $w.$ Let $\rho :=-w^{-1}\Delta _{f}w>0.$ Then

\begin{equation*}
\Delta _{f}w=-\rho w.
\end{equation*}%
This implies that a weighted Poincare inequality on $M$ of the form

\begin{equation}
\int_{M}\rho \phi ^{2}e^{-f}\leq \int_{M}\left\vert \nabla \phi \right\vert
^{2}e^{-f}  \label{s16}
\end{equation}%
is valid for all $\phi $ with compact support in $M.$ Let $\eta >0$ be a
cut-off function so that $\eta =0$ on $B\left( p,r_{0}\right) $ and $\eta =1$
on $M\setminus B\left( p,2r_{0}\right).$ For any function $\phi $ with
compact support in $M$ we have

\begin{equation}
\int_{M}\phi ^{2}e^{-f}\leq 2\int_{M}\left( \phi \eta \right)
^{2}e^{-f}+2\int_{M}\phi ^{2}\left( 1-\eta \right) ^{2}e^{-f}.  \label{s17}
\end{equation}%
Since $\phi \eta $ has support in $M\backslash \Omega $,

\begin{eqnarray*}
\lambda _{1}^{M\backslash \Omega }\left( \Delta _{f}\right) \int_{M}\left(
\phi \eta \right) ^{2}e^{-f} &\leq &\int_{M}\left\vert \nabla \left( \phi
\eta \right) \right\vert ^{2}e^{-f} \\
&\leq &2\int_{M}\left\vert \nabla \phi \right\vert
^{2}e^{-f}+2\int_{M}\left\vert \nabla \eta \right\vert ^{2}\phi ^{2}e^{-f}.
\end{eqnarray*}%
Together with (\ref{s17}), it follows that there exists a constant $\alpha
>0,$ depending on $r_{0}$ and $\lambda _{1}^{M\backslash \Omega }\left(
\Delta _{f}\right) ,$ such that

\begin{equation}
\alpha \int_{M}\phi ^{2}e^{-f}\leq \int_{M}\left\vert \nabla \phi
\right\vert ^{2}e^{-f}+\int_{B\left( p,2r_{0}\right) }\phi ^{2}e^{-f}.
\label{s18}
\end{equation}%
As the weight $\rho >0$ on $M,$ there exists $\beta >0$ so that $\rho \geq
\beta $ on $B\left( p,2r_{0}\right) .$ Therefore, from (\ref{s16}), we get 
\begin{eqnarray*}
\beta \int_{B\left( p,2r_{0}\right) }\phi ^{2}e^{-f} &\leq &\int_{M}\rho
\phi ^{2}e^{-f} \\
&\leq &\int_{M}\left\vert \nabla \phi \right\vert ^{2}e^{-f}.
\end{eqnarray*}%
Combining with (\ref{s18}), we conclude that 
\begin{equation*}
\alpha \left( 1+\beta ^{-1}\right) ^{-1}\int_{M}\phi ^{2}e^{-f}\leq
\int_{M}\left\vert \nabla \phi \right\vert ^{2}e^{-f}
\end{equation*}%
for all $\phi $ with compact support in $M.$ This proves the result.
\end{proof}

\section{Applications to steady Ricci solitons}

In this section, we discuss some applications to steady gradient Ricci
solitons. Recall that complete manifold $\left( M,g\right)$ is a steady
gradient Ricci soliton if there exists a smooth potential function $f$ such
that

\begin{equation*}
\mathrm{Ric}+\mathrm{Hess}\left( f\right) =0.
\end{equation*}%
Hamilton proved that the scalar curvature $S$ satisfies $S+\left\vert \nabla
f\right\vert ^{2}=C$ for some positive constant $C.$ By scaling the metric
if necessary, we may assume that $C=1$ and

\begin{equation}
S+\left\vert \nabla f\right\vert ^{2}=1.  \label{m1}
\end{equation}%
It is known \cite{C} that $S>0$ on $M$ unless the soliton is flat. In
particular, $\left\vert \nabla f\right\vert \leq 1$ on $M$ and $f$ is of
linear growth. The identity (\ref{m1}) together with $S+\Delta f=0$ implies
that

\begin{equation*}
\Delta _{f}\left( f\right) =-1.
\end{equation*}%
Therefore,

\begin{eqnarray*}
\Delta _{f}e^{\frac{f}{2}} &=&\left( \frac{1}{2}\Delta _{f}\left( f\right) +%
\frac{1}{4}\left\vert \nabla f\right\vert ^{2}\right) e^{\frac{f}{2}} \\
&\leq &-\frac{1}{4}e^{\frac{f}{2}}.
\end{eqnarray*}%
It is known (see \cite{LW1}) that the existence of a positive function $u>0$
satisfying $\Delta _{f}u\leq -\frac{1}{4}u$ implies $\lambda _{1}\left(
\Delta _{f}\right) \geq \frac{1}{4}.$ Hence, the weighted Laplacian on a
steady gradient Ricci soliton has positive spectrum. In fact, in \cite{MW1}
it was shown that

\begin{equation}
\lambda _{1}\left( \Delta _{f}\right) =\frac{1}{4}.  \label{m2}
\end{equation}%
Let us also recall some curvature identities on steady gradient Ricci
solitons.

\begin{eqnarray*}
\Delta _{f}S &=&-2\left\vert \mathrm{Ric}\right\vert ^{2} \\
\Delta _{f}\mathrm{Rm} &=&\mathrm{Rm}\star \mathrm{Rm}.
\end{eqnarray*}%
The last identity implies a useful inequality 
\begin{equation*}
\Delta _{f}\left\vert \mathrm{Rm}\right\vert \geq -c\left\vert \mathrm{Rm}%
\right\vert ^{2}.
\end{equation*}%
As before, we fix $p\in M$ and denote 
\begin{equation*}
r\left( x\right) :=d\left( p,x\right).
\end{equation*}%
Also, denote by $C_{0}$, $C$, $c$ constants depending only on the
dimension $n$ and possibly the geometry of $B\left( p,2\right).$

We now estimate the volume of unit balls in $\left( M,g\right).$

\begin{lemma}
\label{Decay} Let $\left( M,g,f\right) $ be a steady gradient Ricci soliton.
Then there exists a constant $C>0$ such that

\begin{equation*}
\mathrm{V}\left( x,1\right) \geq C^{-1}\,e^{-C\sqrt{r\left( x\right) \ln
\left( 1+r\left( x\right) \right) }}\mathrm{V}\left( p,1\right)
\end{equation*}%
for all $x\in M.$
\end{lemma}

\begin{proof}
The proof is inspired by \cite{MW3}, where a similar estimate was given for
shrinking Ricci solitons. Let us denote the volume form in geodesic
coordinates centered at $x$ by 
\begin{equation*}
dV|_{\exp _{x}\left( r\xi \right) }=J\left( x,r,\xi \right) drd\xi
\end{equation*}%
for $r>0$ and $\xi \in S_{x}M,$ the unit tangent sphere at $x.$ Let $%
R:=d\left( p,x\right).$ Without loss of generality, we may assume that $%
R\geq 2.$

Let $\gamma \left( s\right) $ be a minimizing normal geodesic with $\gamma
\left( 0\right) =x$ and $\gamma \left( T\right) \in B\left( p,1\right) $ for
some $T>0.$ By the triangle inequality, we know that 
\begin{equation}
R-1\leq T\leq R+1.  \label{m3}
\end{equation}%
Along $\gamma,$ according to the Laplace comparison theorem,%
\begin{equation*}
m^{\prime }\left( r\right) +\frac{1}{n-1}m^{2}\left( r\right) \leq f^{\prime
\prime }\left( r\right),
\end{equation*}%
where $m\left( r\right) :=\frac{d}{dr}\ln J\left( x,r,\xi \right).$

For arbitrary $k\geq 2,$ multiplying this by $r^{k}$ and integrating from $%
r=0$ to $r=t,$ we have

\begin{equation}
\int_{0}^{t}m^{\prime }\left( r\right) r^{k}dr+\frac{1}{n-1}%
\int_{0}^{t}m^{2}\left( r\right) r^{k}dr\leq \int_{0}^{t}f^{\prime \prime
}\left( r\right) r^{k}dr.  \label{m4}
\end{equation}%
After integrating the first term in (\ref{m4}) by parts and rearranging
terms, we get%
\begin{eqnarray*}
&&m\left( t\right) t^{k}+\frac{1}{n-1}\int_{0}^{t}\left( m\left( r\right) r^{%
\frac{k}{2}}-\left( n-1\right) \frac{k}{2}r^{\frac{k}{2}-1}\right) ^{2}dr \\
&\leq &\frac{\left( n-1\right) k^{2}}{4\left( k-1\right) }%
t^{k-1}+\int_{0}^{t}f^{\prime \prime }\left( r\right) r^{k}dr.
\end{eqnarray*}%
In particular, 
\begin{equation*}
m\left( t\right) \leq \frac{\left( n-1\right) k^{2}}{4\left( k-1\right) }%
\frac{1}{t}+\frac{1}{t^{k}}\int_{0}^{t}f^{\prime \prime }\left( r\right)
r^{k}dr.
\end{equation*}%
Integrating this from $t=1$ to $t=T,$ we obtain for some constant $c$
depending only on $n,$

\begin{equation}
\ln \frac{J\left( x,T,\xi \right) }{J\left( x,1,\xi \right) }\leq c\,k\,\ln
T+A,  \label{m5}
\end{equation}%
where the term $A$ is given by 
\begin{equation*}
A:=\int_{1}^{T}\frac{1}{t^{k}}\left( \int_{0}^{t}f^{\prime \prime }\left(
r\right) r^{k}dr\right) dt.
\end{equation*}

We now estimate $A$ in the right side of (\ref{m5}). Integrating by parts
implies%
\begin{eqnarray*}
A &=&f\left( T\right) -f\left( 1\right) -k\int_{1}^{T}\frac{1}{t^{k}}%
\int_{0}^{t}f^{\prime }\left( r\right) r^{k-1}drdt \\
&=&-\frac{1}{k-1}\left( f\left( T\right) -f\left( 1\right) \right) \\
&&+\frac{k}{k-1}\frac{1}{t^{k-1}}\left( \int_{0}^{t}f^{\prime }\left(
r\right) r^{k-1}dr\right) |_{t=1}^{t=T}.
\end{eqnarray*}%
As $\left\vert f\left( 1\right) \right\vert \leq T+c$ and $\left\vert
f^{\prime }\left( r\right) \right\vert \leq 1,$ it follows that

\begin{equation*}
A\leq \frac{cT}{K}.
\end{equation*}%
We now choose 
\begin{equation*}
k:=\sqrt{\frac{T}{\ln T}}.
\end{equation*}%
It follows from (\ref{m5}) that 
\begin{equation*}
\ln \frac{J\left( x,T,\xi \right) }{J\left( x,1,\xi \right) }\leq ck\ln T+%
\frac{cT}{k}\leq c\sqrt{T\ln T}.
\end{equation*}%
We have thus proved that 
\begin{equation*}
J\left( x,1,\xi \right) \geq e^{-c\sqrt{R\ln R}}J\left( x,T,\xi \right) .
\end{equation*}%
By integrating this over a subset of $S_{x}M$ consisting of all unit tangent
vectors $\xi $ so that $\exp _{x}\left( T\xi \right) \in B\left( p,1\right) $
for some $T,$ it follows that 
\begin{equation*}
\mathrm{A}\left( \partial B\left( x,1\right) \right) \geq e^{-c\sqrt{R\ln R}}%
\mathrm{V}\left( p,1\right) ,
\end{equation*}%
where $R=d\left( p,x\right) \geq 2.$ Clearly, for $\frac{1}{2}\leq t\leq 1,$
a similar estimate holds for $\mathrm{A}\left( \partial B\left( x,t\right)
\right) .$ Therefore, 
\begin{equation*}
\mathrm{V}\left( x,1\right) \geq e^{-c\sqrt{R\ln R}}\mathrm{V}\left(
p,1\right) .
\end{equation*}%
This proves the result.
\end{proof}

From now on, we assume in addition that the potential $f$ is bounded above
by a constant. By adding a constant if necessary, we may assume without loss
of generality that 
\begin{equation}
f\leq 0\text{ \ on }M.  \label{m6}
\end{equation}
Following \cite{CLY}, we now establish a sharp lower bound for the scalar
curvature.

\begin{lemma}
\label{lb}Let $\left( M,g,f\right) $ be a complete steady gradient Ricci
soliton satisfying (\ref{m6}). Then there exists $c>0$ so that 
\begin{equation*}
S\geq c\,e^{f}\text{ on }M.
\end{equation*}
\end{lemma}

\begin{proof}
Using that $\Delta _{f}e^{f}=-Se^{f}$ we compute for $a>0,$%
\begin{eqnarray*}
\Delta _{f}\left( S-ae^{f}\right) &=&-2\left\vert \mathrm{Ric}\right\vert
^{2}+aSe^{f} \\
&\leq &-\frac{2}{n}S^{2}+aSe^{f} \\
&\leq &\frac{na^{2}}{8}e^{2f}.
\end{eqnarray*}%
Note that

\begin{eqnarray*}
\Delta _{f}e^{2f} &=&\left( 2\Delta _{f}\left( f\right) +4\left\vert \nabla
f\right\vert ^{2}\right) e^{2f} \\
&=&2\left( 1-2S\right) e^{2f}.
\end{eqnarray*}%
It follows that for $b>0,$%
\begin{equation*}
\Delta _{f}\left( S-ae^{f}-be^{2f}\right) \leq \left( \frac{na^{2}}{8}%
-2b+4bS\right) e^{2f}.
\end{equation*}%
Let $\phi $ be a smooth cut-off function so that $\phi =1$ on $B\left(
p,R\right) $ and $\phi =0$ on $M\backslash B\left( p,2R\right) .$ We may
assume that 
\begin{eqnarray*}
-\frac{c}{R} &\leq &\phi ^{\prime }\leq 0 \\
\left\vert \phi ^{\prime \prime }\right\vert &\leq &\frac{c}{R^{2}}.
\end{eqnarray*}%
Let 
\begin{equation}
G:=\left( S-ae^{f}-be^{2f}\right) \phi ^{2}.  \label{m7}
\end{equation}%
If $G$ achieves its minimum at $x_{0}$ with $G\left( x_{0}\right) <0,$ then $%
x_{0}\in B\left( p,2R\right) .$ By the Laplacian comparison theorem \cite{WW}
we have on the support of $\phi $ that 
\begin{eqnarray*}
\Delta _{f}r\left( x\right) &\leq &\frac{n-1}{r\left( x\right) }+1 \\
&\leq &c.
\end{eqnarray*}%
It follows that at $x_{0},$%
\begin{eqnarray*}
0 &\leq &\Delta _{f}G \\
&\leq &\left( \frac{na^{2}}{8}-2b+4bS\right) e^{2f}\phi ^{2}-6\left\vert
\nabla \phi \right\vert ^{2}\phi ^{-2}G \\
&&+2\phi ^{-1}\left( \Delta _{f}\phi \right) G+2\left\langle \nabla G,\nabla
\ln \phi \right\rangle \\
&\leq &\left( \frac{na^{2}}{8}-2b+4bS\right) e^{2f}\phi ^{2}+\frac{c}{R}.
\end{eqnarray*}%
Since $G\left( x_{0}\right) <0$ and $f\leq 0,$ we have 
\begin{eqnarray*}
S\left( x_{0}\right) &<&ae^{f\left( x_{0}\right) }+be^{2f\left( x_{0}\right)
} \\
&\leq &a+b.
\end{eqnarray*}%
Hence, we get that 
\begin{equation*}
0\leq \left( \frac{na^{2}}{8}-2b+4b\left( b+a\right) \right) e^{2f\left(
x_{0}\right) }\phi ^{2}+\frac{c}{R}.
\end{equation*}%
Now let $b=\frac{1-2a}{4}$ with $a>0$ sufficiently small. Then the above
inequality yields

\begin{equation}
e^{2f}\phi ^{2}\leq \frac{c}{R}  \label{m8}
\end{equation}%
for some $c>0.$ We claim that

\begin{equation}
G\left( x_{0}\right) \geq -cR^{-\frac{1}{4}}.  \label{m8'}
\end{equation}%
Indeed, if $e^{2f\left( x_{0}\right) }\leq R^{-\frac{1}{2}}$, then (\ref{m7}%
) implies $G\left( x_{0}\right) \geq -cR^{-\frac{1}{4}}$ as claimed. On the
other hand, if $e^{2f\left( x_{0}\right) }>R^{-\frac{1}{2}},$ then (\ref{m8}%
) implies $\phi ^{2}\left( x_{0}\right) \leq c\,R^{-\frac{1}{2}}.$ So from (%
\ref{m7}), $G\left( x_{0}\right) \geq -c\,R^{-\frac{1}{2}}.$ In either case,
(\ref{m8'}) is proved. Certainly, this is true as well if $G\left(
x_{0}\right) \geq 0.$ In conclusion, we have proved that 
\begin{equation*}
G\geq -cR^{-\frac{1}{4}}\text{ on }M.
\end{equation*}%
As $\phi =1$ on $B\left( p,R\right) ,$ one has%
\begin{eqnarray*}
\inf_{B\left( p,R\right) }\left( S-ae^{f}-be^{2f}\right) &\geq &\inf_{M}G \\
&\geq &-cR^{-\frac{1}{4}}.
\end{eqnarray*}%
Letting $R\rightarrow \infty ,$ we conclude that $S-ae^{f}-be^{2f}\geq 0$ on 
$M.$ This proves the result.
\end{proof}

We now prove the main result of this section.

\begin{theorem}
\label{Rm_bd}Let $\left( M,g,f\right) $ be a complete steady gradient Ricci
soliton satisfying (\ref{m6}). If%
\begin{equation*}
\lim_{x\rightarrow \infty }\left\vert \mathrm{Rm}\right\vert \left( x\right)
r\left( x\right) =0,
\end{equation*}%
then there exists $c>0$ such that 
\begin{equation*}
\left\vert \mathrm{Rm}\right\vert \left( x\right) \leq c\left( 1+r\left(
x\right) \right) ^{3\left( n+1\right) }e^{-r\left( x\right) }\text{ \ on }M.
\end{equation*}
\end{theorem}

\begin{proof}
Since $S\leq c\left( 1+r\left( x\right) \right) ^{-1},$ by Lemma \ref{lb} $f$
is proper and

\begin{equation*}
f\left( x\right) \leq -c_{1}\ln \left( 1+r\left( x\right) \right) +c_{2}.
\end{equation*}
For $\sigma >8$ fixed and to be specified later, we define 
\begin{equation*}
F:=f+2\ln \left( -f+\sigma \right) .
\end{equation*}%
One checks directly that 
\begin{equation*}
\Delta _{F}\left( f\right) =-1+2\left\vert \nabla f\right\vert ^{2}\left(
-f+\sigma \right) ^{-1}.
\end{equation*}%
Hence,

\begin{eqnarray*}
\Delta _{F}e^{\frac{f}{2}} &=&\frac{1}{2}\left( -1+2\left\vert \nabla
f\right\vert ^{2}\left( -f+\sigma \right) ^{-1}+\frac{1}{2}\left\vert \nabla
f\right\vert ^{2}\right) e^{\frac{f}{2}} \\
&\leq &-\rho e^{\frac{f}{2}},
\end{eqnarray*}%
where 
\begin{eqnarray}
\rho &=&\frac{1}{4}-\left( -f+\sigma \right) ^{-1}  \label{rho} \\
&\geq &\frac{1}{4}-\frac{1}{\sigma }.  \notag
\end{eqnarray}%
It is well known (see \cite{LW}) that this implies an estimate for the
bottom of spectrum of the weighted Laplacian $\Delta _{F}:=\Delta
-\left\langle \nabla F,\nabla \right\rangle $ of the form

\begin{equation*}
\lambda _{1}\left( \Delta _{F}\right) \geq \frac{1}{4}-\frac{1}{\sigma }>0.
\end{equation*}%
Applying (\ref{LW}) to the Green's function $\bar{G}\left( x,y\right) $ of $%
\Delta _{F}$ we get%
\begin{equation}
\int_{B\left( p,R+1\right) \backslash B\left( p,R\right) }\bar{G}^{2}\left(
p,y\right) e^{-F\left( y\right) }dy\leq Ce^{-2\sqrt{\lambda _{1}\left(
\Delta _{F}\right) }R}  \label{m9}
\end{equation}%
for any $R\geq 1$. For any $x\in M\backslash B\left( p,2\right) $ we have by
the triangle inequality that 
\begin{equation*}
B\left( x,1\right) \subset B\left( p,r\left( x\right) +2\right) \backslash
B\left( p,r\left( x\right) -1\right) .
\end{equation*}%
It is easy to see from (\ref{m9}) that 
\begin{equation}
\int_{B\left( x,1\right) }\bar{G}^{2}\left( p,y\right) e^{-F\left( y\right)
}dy\leq Ce^{-\left( 1-\frac{4}{\sigma }\right) r\left( x\right) }
\label{m10}
\end{equation}%
for any $x\in M\backslash B\left( p,2\right) $. By (\ref{m1}) we get 
\begin{eqnarray*}
\left\vert \nabla F\right\vert &=&\left( 1-\frac{2}{\left( -f+\sigma \right) 
}\right) \left\vert \nabla f\right\vert \\
&\leq &1.
\end{eqnarray*}%
We note that the Bakry-Emery tensor associated to the weight $F$ is 
\begin{eqnarray*}
\mathrm{Ric}_{F} &=&\mathrm{Ric}_{f}+\mathrm{Hess}\left( 2\ln \left(
-f+\sigma \right) \right) \\
&=&-2\left( -f+\sigma \right) ^{-1}\mathrm{Hess}\left( f\right) -2\left(
-f+\sigma \right) ^{2}\nabla f\otimes \nabla f.
\end{eqnarray*}%
As $M$ has bounded Ricci curvature $\left\vert \mathrm{Ric}\right\vert \leq
c,$ it is clear from above that 
\begin{equation*}
\mathrm{Ric}_{F}\geq -\left( n-1\right) \bar{K}
\end{equation*}%
for some $\bar{K}>0$ independent of $\sigma.$ Since $\Delta _{F}\bar{G}%
\left( p,\cdot \right) =0$ on $B\left( x,2\right),$ by (\ref{b2}) we get a
gradient estimate of the form 
\begin{equation*}
\sup_{y\in B\left( x,1\right) }\left\vert \nabla \ln \bar{G}\right\vert
\left( p,y\right) \leq c
\end{equation*}%
for a constant $c>0$ independent of $\sigma.$ Integrating this
estimate, one sees that 
\begin{equation*}
\bar{G}\left( p,x\right) \leq c\bar{G}\left( p,y\right)
\end{equation*}%
for all $y\in B\left( x,1\right).$ Hence, (\ref{m10}) implies that%
\begin{equation*}
\bar{G}\left( p,x\right) \leq C\mathrm{V}^{-\frac{1}{2}}\left( x,1\right)
e^{-\left( \frac{1}{2}-\frac{2}{\sigma }\right) r\left( x\right) +\frac{1}{2}%
F\left( x\right) }.
\end{equation*}
Since $F=f+2\ln \left( -f+\sigma \right),$ we obtain 
\begin{equation}
\bar{G}\left( p,x\right) \leq C\mathrm{V}^{-\frac{1}{2}}\left( x,1\right)
e^{-\left( \frac{1}{2}-\frac{2}{\sigma }\right) r\left( x\right) +\frac{1}{2}%
f\left( x\right) }\left( -f\left( x\right) +\sigma \right)  \label{m12}
\end{equation}%
for all $x\in M\backslash B\left( p,2\right).$

On the other hand, recall that $\Delta _{f}\left\vert \mathrm{Rm}\right\vert
\geq -c\left\vert \mathrm{Rm}\right\vert ^{2}$and $\Delta _{f}\left(
-f\right) =1$. It follows that 
\begin{eqnarray*}
\Delta _{f}\left( \left\vert \mathrm{Rm}\right\vert \left( -f+\sigma \right)
\right) &\geq &-c\left\vert \mathrm{Rm}\right\vert ^{2}\left( -f+\sigma
\right) +\left\vert \mathrm{Rm}\right\vert \\
&&+2\left\langle \nabla \left\vert \mathrm{Rm}\right\vert ,\nabla \left(
-f+\sigma \right) \right\rangle \\
&=&\left\vert \mathrm{Rm}\right\vert \left( 1-c\left\vert \mathrm{Rm}%
\right\vert \left( -f+\sigma \right) \right) \\
&&+2\left\langle \nabla \left( \left\vert \mathrm{Rm}\right\vert \left(
-f\right) \right) ,\nabla \ln \left( -f+\sigma \right) \right\rangle \\
&&-2\left\vert \mathrm{Rm}\right\vert \left( -f+\sigma \right)
^{-1}\left\vert \nabla f\right\vert ^{2}.
\end{eqnarray*}%
So the function $w:=\left\vert \mathrm{Rm}\right\vert \left( -f+\sigma
\right) $ satisfies 
\begin{equation*}
\Delta _{F}w\geq \left\vert \mathrm{Rm}\right\vert \left( \frac{1}{2}%
-c\left\vert \mathrm{Rm}\right\vert \left( -f+\sigma \right) \right).
\end{equation*}%
As $\left\vert \mathrm{Rm}\right\vert \left( x\right) r\left( x\right)
=o\left( 1\right),$ there exists $c>0$ independent of $\sigma $ so that 
\begin{equation*}
\Delta _{F}w\geq 0\text{ on }M\backslash B\left( p,c\sigma \right).
\end{equation*}%
Fix $A>0$ large enough with
\begin{equation*}
w\left( x\right) -A\bar{G}\left( p,x\right) <0\text{ \ for all }x\in
\partial B\left( p,c\sigma \right).
\end{equation*}%
By (\ref{m6}), (\ref{m12}) and Lemma \ref{Decay} we get that $\bar{G}\left(
p,x\right) $ converges to zero as $x\rightarrow \infty.$ Then the function $%
w-A\bar{G}\left( p,\cdot \right) $ is $F$-subharmonic on $M\backslash
B\left( p,c\sigma \right),$ converges to zero at infinity and is
non-positive on $\partial B\left( p,c\sigma \right).$ The maximum
principle implies that 
\begin{equation*}
w\left( x\right) -A\bar{G}\left( p,x\right) <0\text{ \ for all \ }x\in
M\backslash B\left( p,c\sigma \right).
\end{equation*}%
Combining with (\ref{m12}) we get 
\begin{equation*}
\left\vert \mathrm{Rm}\right\vert \left( x\right) \leq C(\sigma)\mathrm{V}^{-%
\frac{1}{2}}\left( x,1\right) e^{-\left( \frac{1}{2}-\frac{2}{\sigma }%
\right) r\left( x\right) +\frac{1}{2}f\left( x\right) }.
\end{equation*}%
Lemma \ref{lb} implies $e^{f\left( x\right) }\leq \left\vert \mathrm{Rm}%
\right\vert \left( x\right).$ Hence we get from above that%
\begin{equation*}
\left\vert \mathrm{Rm}\right\vert \left( x\right) \leq C\left( \sigma
\right) \mathrm{V}^{-1}\left( x,1\right) e^{-\left( 1-\frac{4}{\sigma }%
\right) r\left( x\right) }.
\end{equation*}
Together with Lemma \ref{Decay} this proves that for given $\varepsilon >0$
there exists $C\left( \varepsilon \right) >0$ so that 
\begin{equation}
\left\vert \mathrm{Rm}\right\vert \left( x\right) \leq C\left( \varepsilon
\right) e^{-\left( 1-\varepsilon \right) r\left( x\right) }.  \label{m13}
\end{equation}%
To finish the proof of the theorem, we use Theorem \ref{Poisson1}. Note that
by (\ref{m13}) we have $\left\vert \mathrm{Rm}\right\vert \left( x\right)
\leq ce^{-\frac{3}{4}r\left( x\right) }$ for all $x\in M.$ Hence, solving
the Poisson equation $\Delta _{f}u=-c\left\vert \mathrm{Rm}\right\vert ^{2}$
by choosing $\alpha =1$ in Theorem \ref{Poisson1} and noticing that
$K=0$ and $a=b=1$ due to our normalization, 
we obtain a solution $u$ such that

\begin{gather}
\left\vert u\right\vert \left( x\right) \leq C\left( \int_{2r\left( x\right)
}^{\infty }\omega \left( t\right) dt+r\left( x\right) \omega \left( r\left(
x\right) \right) \right)  \label{m14} \\
+C\left( 1+r\left( x\right) \right) ^{n+1}e^{-\frac{1}{2}r\left( x\right) +%
\frac{1}{2}f\left( x\right) }\mathrm{V}^{-\frac{1}{2}}\left( x,1\right)
\int_{0}^{r\left( x\right) }\omega \left( t\right) e^{t}dt,  \notag
\end{gather}%
where 
\begin{equation*}
\omega \left( t\right) =ce^{-\frac{3}{2}t}.
\end{equation*}
By (\ref{m13}) and standard comparison geometry we know that 
\begin{equation*}
\mathrm{V}\left( x,1\right) \geq c^{-1}\left( 1+r\left( x\right) \right)
^{-n-1}.
\end{equation*}%
Therefore (\ref{m14}) implies that 
\begin{eqnarray}
\left\vert u\right\vert \left( x\right) &\leq &C\left( 1+r\left( x\right)
\right) ^{\frac{3}{2}\left( n+1\right) }e^{-\frac{1}{2}r\left( x\right) +%
\frac{1}{2}f\left( x\right) }  \label{m15} \\
&\leq &C\left( 1+r\left( x\right) \right) ^{\frac{3}{2}\left( n+1\right)
}e^{-\frac{1}{2}r\left( x\right) }\sqrt{\left\vert \mathrm{Rm}\right\vert
\left( x\right) },  \notag
\end{eqnarray}%
where in the last line we have used Lemma \ref{lb}. Since $\Delta
_{f}\left\vert \mathrm{Rm}\right\vert \geq -c\left\vert \mathrm{Rm}%
\right\vert ^{2}$ and $\Delta _{f}u=-c\left\vert \mathrm{Rm}\right\vert ^{2},$ 
by the maximum principle $\left\vert \mathrm{Rm}\right\vert \leq u$ on $M.$
In conclusion, (\ref{m15}) implies that%
\begin{equation*}
\left\vert \mathrm{Rm}\right\vert \left( x\right) \leq C\,\left( 1+r\left(
x\right) \right) ^{3\left( n+1\right) }e^{-r\left( x\right) }.
\end{equation*}%
This proves the theorem.
\end{proof}

Finally we point out that stronger results can be obtained by assuming the
sectional curvature is non-negative. First, we recall a result from \cite{D}
and \cite{CMM}. For completeness, a simple proof is provided here.

\begin{proposition}
\label{Classif}Let $\left( M^n,g,f\right) $ be an $n$-dimensional complete
non-flat steady gradient Ricci soliton with non-negative sectional
curvature. Assume that the scalar curvature is integrable on $M.$ Then $%
\left( M,g\right) $ is isometric to a quotient of $\mathbb{R}^{n-2}\times
\Sigma,$ where $\Sigma $ is the cigar soliton.
\end{proposition}

\begin{proof}
We may assume $n\geq 3$ as otherwise the result is known. We note that $%
2\left\vert \mathrm{Ric}\right\vert ^{2}\leq S^{2}.$ Indeed, the fact that $%
M $ has nonnegative sectional curvature implies that the eigenvalues $%
\lambda_{i}$ of the Ricci curvature satisfy $\sum_{j\neq i}\lambda _{j}\geq
\lambda_{i}.$ So $S\geq 2\lambda _{i}$ and

\begin{equation}
2\left\vert \mathrm{Ric}\right\vert ^{2}=2\sum_{i}\lambda _{i}^{2}\leq
\sum_{i}\left( \lambda _{i}S\right) =S^{2}.  \label{m30}
\end{equation}%
Hence, using the cut-off function $\phi :=\left( \frac{R-d\left( p,x\right) 
}{R}\right) _{+}$ with support in $B\left( p,R\right) $ we get 
\begin{eqnarray*}
0 &\leq &\int_{M}\left( S^{2}-2\left\vert \mathrm{Ric}\right\vert
^{2}\right) \phi ^{2} \\
&=&-\int_{M}S\left( \Delta f\right) \phi ^{2}+2\int_{M}R_{ij}f_{ij}\phi ^{2}
\\
&=&\int_{M}\left\langle \nabla S,\nabla f\right\rangle \phi
^{2}-2\int_{M}\left( \nabla _{j}R_{ij}\right) f_{i}\phi ^{2} \\
&&+\int_{M}\left\langle \nabla f,\nabla \phi ^{2}\right\rangle
S-2\int_{M}R_{ij}f_{i}\left( \phi ^{2}\right) _{j} \\
&\leq &\frac{c}{R}\int_{B\left( p,R\right) }S,
\end{eqnarray*}%
where in the last line we have used that $2\nabla _{j}R_{ij}=\nabla _{i}S$
by the Bianchi identity. Therefore, by letting $R\rightarrow \infty,$ we
conclude $2\left\vert \mathrm{Ric}\right\vert ^{2}=S^{2}$ on $M.$ In
particular, from (\ref{m30}) we see that either $\lambda _{i}=0$ or $\lambda
_{i}=\frac{1}{2}S.$ If $\lambda _{i}=\frac{1}{2}S$ for all $i$ at all
points, then $M$ is Einstein and flat. So $\lambda _{i}=0$ for some $i$ at
some point. Applying Hamilton's strong maximum principle, the result follows.
\end{proof}

Combining the proposition with Theorem \ref{Rm_bd} we get the following.

\begin{corollary}
Let $\left( M,g,f\right) $ be a complete non-flat steady gradient Ricci
soliton with non-negative sectional curvature. Assume that the scalar
curvature decays faster than linear, that is, 
\begin{equation*}
\lim_{x\rightarrow \infty } r\left( x\right) S\left( x\right)=0.
\end{equation*}%
Then $\left( M,g\right) $ is isometric to a quotient of $\mathbb{R}%
^{n-2}\times \Sigma,$ where $\Sigma $ is the cigar soliton.
\end{corollary}

\begin{proof}
By \cite{CN}, the function $f$ satisfies

\begin{equation*}
-r\left( x\right) +c_{1}\leq f\left( x\right) \leq -cr\left( x\right) +c_{1}.
\end{equation*}%
Theorem \ref{Rm_bd} implies that the curvature decays exponentially, that
is, 
\begin{equation*}
\left\vert \mathrm{Rm}\right\vert \left( x\right) \leq c\left( 1+r\left(
x\right) \right) ^{3\left( n+1\right) }e^{-r\left( x\right) }.
\end{equation*}%
This shows that $S\in L^{1}\left( M\right) $ and the result follows from
Proposition \ref{Classif}.
\end{proof}

\end{document}